\documentclass[11pt]{article}
\author{Bing-Long Chen\footnote {The research
partially supported by Grants 2005-34000-3171404 and
2006-34000-1131040. \newline AMS Mathematics Subject
Classification
Numbers: Primary 53c44; secondary 35k55. }  and   Le Yin \\[8pt]}
\title{\textbf{Uniqueness and Pseudolocality Theorems of the Mean Curvature Flow}}
\date{Revised at June 5, 2007}
\usepackage{latexsym}
\usepackage{amsmath}
\usepackage{amssymb,amscd}
\usepackage[mathscr]{eucal}
\newtheorem{thm}{Theorem}[section]
\newtheorem{cor}[thm]{Corollary}
\newtheorem{lem}[thm]{Lemma}
\newtheorem{prop}[thm]{Proposition}

\numberwithin{equation}{section}
\newenvironment{pf}{{\noindent \it  Proof.}}{{\hfill$\Box$}\\}
\begin{document}
\maketitle

\begin{abstract}

Mean curvature flow evolves isometrically immersed base manifolds
$M$ in the direction of their mean curvatures in an ambient
manifold $\bar{M}$. If the base manifold $M$ is compact,
  the short time existence and uniqueness of the mean curvature flow  are well-known. For  complete isometrically
  immersed submanifolds of arbitrary codimensions, the existence
  and uniqueness   are still unsettled even in the Euclidean space. In this paper,
  we solve the uniqueness problem affirmatively
  for the mean curvature flow of general codimensions and general ambient manifolds. In the second
  part of the paper, inspired by the Ricci flow, we prove a pseudolocality theorem of mean curvature flow.
  As a
  consequence, we obtain a strong uniqueness theorem, which removes
  the assumption on the boundedness  of the second fundamental form of the solution.
\end{abstract}

  \section{Introduction}
\vskip 0.5cm

\qquad  Let $(\bar{M}^{\bar{n}}, \bar{g})$ be a complete
Riemannian (compact or noncompact) manifold,  and $X_0:
(M^{n},g)\rightarrow \bar{M}^{\bar{n}}$ be an isometrically
immersed Riemannian manifold. For any fixed point $x_0\in M^{n}$,
$X,Y \in T_{x_0}M^{n}$, the second fundamental form $II$ at $x_0$
is defined by
$II(X,Y)=\bar{\nabla}_{\tilde{X}}\tilde{Y}-\nabla_{\tilde{X}}\tilde{Y}=
(\bar{\nabla}_{\tilde{X}}\tilde{Y})^{\bot}$, where  $M^{n}$ is
regarded as a submanifold of $\bar{M}$ locally by the isometry
$X_0$,  $\bar{\nabla}$ and $\nabla$ are the covariant derivatives
of $\bar{g}$ and $g$ respectively, $\tilde{X},\tilde{Y}$ are any
smooth extensions of $X$ and $Y$ on $\bar{M}^{\bar{n}}$. In local
coordinate system $\{x^{1},x^{2},\cdots,x^{n}\}$ on $M^{n}$,
denote the second fundamental form by
$h_{ij}=II(\frac{\partial}{\partial
x^{i}},\frac{\partial}{\partial x^{j}})$ and the mean curvature by
 $H=g^{ij}h_{ij}$.
 The mean curvature flow (MCF) is a deformation
  $X_t:M^{n}\rightarrow \bar{M}^{\bar{n}}$  of $X_0$ in the direction of the mean curvature $H$

$$
\frac{\partial}{\partial t}X(x,t)=H(x,t), \ \ \ \ \ \ \ \ \
\text{for}\ x\in M^{n} \text{ and }\  t\geq 0, \eqno{(1.1)}
$$
with $X(x,0)=X_0(x)$, where $M^{n}$ is equipped with the induced
metric from  $X(\cdot,t):M^{n}\rightarrow \bar{M}^{\bar{n}}$ and
$H(x,t)$ is the corresponding  mean curvature.  We can write (1.1)
in another form
$$
\frac{\partial}{\partial t}X(x,t)=\triangle X(x,t), \ \ \ \ \ \ \
\ \ \text{for}\ x\in M^{n} \text{ and }\  t\geq 0, \eqno{(1.2)}
$$ where
$\triangle
X^{\alpha}(x,t)=g^{ij}(x,t)(\frac{\partial^{2}X^{\alpha}}{\partial
    x^{i}\partial
   x^{j}}-\Gamma^{k}_{ij}\frac{\partial X^{\alpha}}{\partial
   x^{k}}+\Gamma^{\alpha}_{\beta\gamma}\frac{\partial X^{\beta}}{\partial x^{i}}
   \frac{\partial X^{\gamma}}{\partial x^{j}}) $ is the harmonic
   map Laplacian from the manifold $(M^{n},g_{ij}(\cdot,t))$ to
   $(\bar{M}^{\bar{n}}, \bar{g})$, and $g_{ij}(\cdot,t)$ is the
   induced metric from the inclusion map $X(\cdot,t)$.

   Various weak solutions to the MCF
   have been studied in the past 30 years by many mathematicians with
   different approaches, e.g. Brakke
   solutions, the level set solutions, etc. The existence, uniqueness and non-uniqueness of
   weak solutions for Euclidean (non)smooth hypersurface have been
   extensively studied. In this paper, motivated by  geometric
   applications, we consider the classical solutions  in  general ambient Riemannian manifolds.

   When
   $M^{n}$ is compact, the MCF (1.1) has a unique
   short time solution, since (1.2) is a (degenerate) quasi-linear parabolic
   equation. For codimensional one complete immersed
   local Lipschitz hypersurfaces  in the Euclidean
  space, we refer the readers to see \cite{EH}.
   For  submanifolds of arbitrary codimensions in a general ambient Riemannian manifold, the short time existence and the
   uniqueness
   of (1.1) have not been  established
   in the literature.
   In
   this paper, we  deal with the  uniqueness problem of
   the mean curvature flow and derive the pseudolocality
   estimate.

    The first main theorem of this paper is the following
   \begin{thm} Let $(\bar{M}^{\bar{n}},\bar{g})$ be a complete  Riemannian manifold
  of dimension $\bar{n}$ such that the curvature and its covariant derivatives up to order 2 are bounded and the injectivity radius
  is bounded from below by a positive constant, i.e. there are
 constants $\bar{C}$ and $\bar{\delta}$ such that
 $$|\bar{R}m|+|\bar{\nabla}\bar{R}m|+|\bar{\nabla}^{2}\bar{R}m|(x)\leq
 \bar{C}, \ \ \ \
 inj(\bar{M}^{\bar{n}},x)>\bar{\delta}>0,  $$
 for all $
 x\in \bar{M}^{\bar{n}}.$
  Let $X_0:M^{n}\rightarrow
 \bar{M}^{\bar{n}}$ be an isometrically immersed Riemannian manifold with bounded
 second fundamental form in $\bar{M}^{\bar{n}}$. Suppose $X_1(x,t)$
 and $X_2(x,t)$ are two solutions to the mean curvature flow (1.1)
 on $M^n\times[0,T]$ with the same $X_0$ as initial data and with
 bounded second fundamental forms on $[0,T]$. Then $X_1(x,t)=
 X_2(x,t)$ for all $(x,t)\in M^{n}\times [0,T]$.
   \end{thm}

 We remark that the uniqueness of the Ricci flow has been
established by Zhu and the first author in \cite{CZ2}. More
precisely, it was proved in \cite{CZ2} that the solutions of the
Ricci flow in the class of bounded curvature with the same initial
data are unique.  We refer the reader to see an interesting
application of  this uniqueness theorem to  the theory of the Ricci
flow with surgery in dimension three  and four\cite{CZ1}. We hope
this MCF uniqueness theorem will also play roles in the theory of
the mean curvature flow with surgery.

 Since the MCF is degenerate in
 tangent directions, it is not a strictly parabolic system.
 In order to apply the standard theory of strict parabolic
 equations, we use the De Turck trick \cite{De}.  The idea is to pull
 back the MCF through a family of diffeomorphisms of the
 base manifold $M^{n}$ generated
 by solving a harmonic map flow coupled with the MCF,
  this gives us the so-called mean curvature De Turck flow, which is a strict
 parabolic system. Then we apply the uniqueness of the strict parabolic
 system. The issue is not quite straight forward as it seems.
 Because before applying the uniqueness theorem of a strict parabolic system on a noncompact manifold, we
   encounter two analytic difficulties. The first one is that we need to establish a short time
   existence for the harmonic map flow between complete manifolds. The second one is to get a priori estimates for the
   harmonic map flow so that after pulling back, the solutions to the strictly parabolic system still satisfy
   suitable smooth or growth
   conditions.

   In the classical theory of the harmonic map flow, people usually would like to impose
    certain convexity conditions to ensure the
   existence (e.g. the negative curvature condition \cite{ES} or convex condition \cite{DL}).
    We observed that in \cite{CZ2} the condition of
   injectivity radius bounded from below by a positive constant ensures certain uniform
   (local) convexity and  this is sufficient to give
   the short time existence and  a priori estimates for the
   harmonic map flow.   Note that the MCF is a
    kind of harmonic map flow with varying base metrics.
    In order to deal with the a priori estimates for MCF  and harmonic map
    flow coupled with MCF, we have to consider the general harmonic map flow.
     These estimates have been dealt with systematically in this
    paper(Sections 2, 3 and 4).

   Note that the injectivity radius of a Riemannian manifold
   with bounded curvature
   may
   decay  exponentially. In the Ricci flow case  \cite{CZ2},
    since we only have the curvature bound,
    we need make more effort to
   overcome this difficulty.

    The difference of Theorem 1.1  with \cite{CZ2} is between the extrinsic and intrinsic
   geometries. In the present case, instead of the metrics as in the Ricci
   flow, we consider
   the equation of the position function.

   As a direct consequence of
Theorem 1.1, we have
   \begin{cor}
   Let $(\bar{M}^{\bar{n}},\bar{g})$ be assumed as in Theorem 1.1 and  $X_t:M^{n}\rightarrow
\bar{M}^{\bar{n}}$
 be a solution to the mean curvature flow (1.1)
on $M^n\times[0,T]$ with bounded second fundamental forms on
$[0,T]$,  and with complete isometric immersed $X_0:M\rightarrow
M$ initial data.  Let $\bar{\sigma}$ be an isometry of
$(\bar{M}^{\bar{n}},\bar{g})$ such that there is an isometry
$\sigma$ of $(M^{n},g)$ to itself satisfying
$$(\bar{\sigma}\circ X_0)(x)=(X_0\circ \sigma)(x) \eqno{(1.3)}
$$
for all $
 x\in {M}^{n}.$ Then  we have
$$(\bar{\sigma}\circ X_t)(x)=(X_t\circ \sigma)(x) \eqno{(1.4)}
$$
for all $(x,t)\in M^{n}\times[0,T].$  In particular, the isometry
subgroup of $(M^{n},g)$  induced by an isometry subgroup of
$(\bar{M}^{\bar{n}},\bar{g})$ at initial time remains to be an
isometry subgroup of $(M^{n},g_t)$ for any $t\in [0,T].$
   \end{cor}

   From the PDE point of view, it is a natural condition in Theorem 1.1 to assume that the
   second fundamental form of the solution is bounded.  In the last
   part of the paper, we try to remove this condition. We remark
   that in \cite{ChZ}, Chou and Zhu  have obtained the strong uniqueness of the curve shortening flow
   for the locally Lipschitz continuous properly embedded curve
   whose two ends are presentable as graphs over semi-infinite
   line. Our strong uniqueness theorem is the following
\begin{thm} Let $\bar{M}$ be an $\bar{n}$-dimensional complete Riemannian manifold
satisfying $\sum\limits_{i=0}^{3}|\bar{\nabla}^{i}\bar{R}m|\leq
c_0^{2}$ and $inj(\bar{M})\geq i_0>0$. Let $X_0:M\rightarrow
 \bar{M}$ be an $n$-dimensional isometrically properly embedded  submanifold with bounded
 second fundamental form in $\bar{M}$. We assume $X_0(M)$ is uniform graphic with some radius $r>0.$  Suppose $X_1(x,t)$
 and $X_2(x,t)$ are two smooth  solutions to the mean curvature flow (1.1)
 on $M\times[0,T_0]$ properly embedded  in $\bar{M} $ with the same $X_0$ as initial data.
  Then there is $0<T_1\leq T_0$ such that $X_1(x,t)=
 X_2(x,t)$ for all $(x,t)\in M\times [0,T_1]$.
\end{thm}

Here roughly speaking, uniform graphic with radius $r>0$ means
that for any $x_0\in X_0(M),$ $X_0(M)\cap B_{\bar{M}}(x_0,r)$ is a
graph.
 We say a
submanifold $M\subset \bar{M}$ is properly embedded in a ball
$B_{\bar{M}}(x_0, r_0)$ if either $M$ is closed or $\partial M$ has
distance $\geq r_0$ from $x_0.$ A submanifold $M\subset \bar{M}$ is
said to be properly embedded in (complete manifold)$\bar{M}$ if
either $M$ is closed or there is an $x_0\in \bar{M}$ such that $M$
is properly embedded in $B_{\bar{M}}(x_0, r_0)$ for any $r_0>0.$

The strong uniqueness theorem was proved as a consequence of
Theorem 1.1 and pseudolocality theorem.

The pseudolocality theorem says that the behavior of the solution
at a point can be controlled by the initial data of  nearby
points,
 no matter the solution or initial data outside the neighborhood
behaviors like. Precisely the following theorem is proved in this
paper:
\begin{thm} Let $\bar{M}$ be an $\bar{n}$-dimensional manifold
satisfying $\sum\limits_{i=0}^{3}|\bar{\nabla}^{i}\bar{R}m|\leq
c_0^{2}$ and $inj(\bar{M})\geq i_0>0$. Then for every $\alpha>0$
there exist $\varepsilon>0$, $\delta>0$ depending only on the
constants $\bar{n}$, $c_0$ and $i_0$ with the following property.
Suppose we have a smooth solution to the mean curvature flow
$M_t\subset\bar{M}$ properly embedded in $B_{\bar{M}}(x_0, r_0)$
for $t\in [0,T]$,  where $0<T\leq \varepsilon^{2} r_0^{2}$, and
assume that at time zero, $M_0$ is a local $\delta$- Lipschitz
graph of radius $r_0$ at $x_0\in {M}$ with $r_0\leq\frac{i_0}{2}$.
Then we have an estimate of the second fundamental form
$$
|A|(x,t)^{2}\leq\frac{\alpha}{t}+(\varepsilon r_0)^{-2}
$$
on $B_{\bar{M}}(x_0,\varepsilon r_0)\cap M_t$, for any $t\in
[0,T]$.
\end{thm}

We refer the reader to see the precise definition of $\delta$-
Lipschitz graph in section 7.  The third covariant derivative of
the curvature is a technical assumption which could be improved,
we assume it only for simplicity. For most of interesting cases,
we have all covariant derivative bounds.

 We remark that for
codimension one uniformly local Lipschitz hypersurface in
Euclidean space, the estimate was firstly derived by Ecker and
Huisken \cite{EH}. For higher codimension case, under an
additional condition which assumes that the submanifold is
compact, the estimate was proved by M.T.Wang\cite{W1}. In
codimension one case \cite{EH}, the constant $\delta$ in Theorem
1.4 does not need to be small;
   however, in higher codimension case, as noted by \cite{W1}, the
smallness assumption is necessary in view of the example of Lawson
and Osserman \cite{LO}.
 The
strategy of the  proofs of \cite{EH} \cite{W1} is to find  a
suitable gradient function.   The philosophy is that this gradient
function will serve as the lower order quantity as in the
Bernstein trick, and the second fundamental form is the higher
order quantity, then apply the maximum principle.

  Our approach
is completely different.  This  approach can be regarded as an
integral version of Bernstein trick. It is a mean curvature flow
analogue of the corresponding estimate in Ricci flow  given by
Perelman \cite{P1}.

As a nontrivial corollary of Theorem 1.4, we have

\begin{cor} Let $\bar{M}$ be an $\bar{n}$-dimensional complete manifold
satisfying $\sum_{i=0}^{3}|\bar{\nabla}^{i}\bar{R}m|\leq c_0^{2}$
and $inj(\bar{M})\geq i_0>0$. Let $X_0:M\rightarrow
 \bar{M}$ be an $n$-dimensional isometrically properly embedded  submanifold with bounded
 second fundamental form $|A|\leq c_0$ in $\bar{M}$. We assume $M_0=X_0(M)$ is uniform graphic with some radius $r>0.$
  Suppose $X(x,t)$
 is a smooth solution to the mean curvature flow (1.1)
 on $M\times[0,T_0]$ properly embedded in  $\bar{M}$ with $X_0$ as initial data.
  Then there is $T_1>0$ depending upon $c_0, i_0, r$ and
  the dimension $\bar{n}$ such that $$|A|(x,t)\leq 2c_0$$ for all $x\in M,$ $0\leq t\leq \min\{T_0,T_1\}$.
\end{cor}

 This paper is
organized as follows. In section 2, we derive the injectivity
radius estimate of an immersed manifold and some preliminary
estimates for a general harmonic map flow. In section 3, the
higher derivative estimates for the MCF are derived.  In Section
4, we study the harmonic map flow coupled with the MCF. In Section
5, we deal with the uniqueness theorem of the mean curvature De
Turck flow. In section 6, we prove the uniqueness Theorem 1.1 and
Corollary 1.2. In section 7, we establish the pseudolocality
theorems 1.4,1.5 and prove the strong uniqueness theorem 1.3.

 \vskip 0.3cm

We are grateful to Professor Xi-Ping Zhu for useful conversations
and encouragement. The second author would like to thank Professor
Nai-Chung Leung and Professor Luen-Fai Tam for their constant
teaching and encouragement, and Professor Mu-Tao Wang for very
helpful discussions.

\vskip 0.8cm
\section{Preliminary estimates}
\qquad In the first part of this section, we will derive the
injectivity radius estimate for isometrically immersed manifold
$M^{n}.$
\begin{thm} Let $(\bar{M}^{\bar{n}},\bar{g})$ be a complete Riemannian manifold
  of dimension $\bar{n}$ with bounded curvature and the injectivity radius
  is bounded from below by a positive constant, i.e. there are
 constants $\bar{C}$ and $\bar{\delta}$ such that
 $$|\bar{R}m|(x)\leq
 \bar{C} \ \ \ \ \mbox{and} \ \ \ \
 inj(\bar{M}^{\bar{n}},x)\geq\bar{\delta}>0, \ \ \ \mbox{for all }
 x\in \bar{M}^{\bar{n}}.\eqno(2.1)$$ Let $X:M^{n}\rightarrow \bar{M}^{\bar{n}}$ be a complete isometrically immersed manifold with bounded
 second fundamental form $|h_{ij}^{\alpha}|\leq C$ in $\bar{M}^{\bar{n}}$, then there is a positive constant
 $\delta=\delta(\bar{C},\bar{\delta},C,\bar{n})$ such that the injectivity radius of $M^{n}$ satisfies
  $$
 inj({M}^{{n}},x)\geq{\delta}>0, \ \ \ \mbox{for all }
 x\in {M}^{{n}}.\eqno(2.2)$$
   \end{thm}
   \begin{pf} Fix $x_0\in M^{n}$, let $\{y^{1},y^{2},\cdots,y^{\bar{n}}\}$ and $\{x^{1},x^{2},\cdots,x^{n}\}$
   be any two local coordinates of $\bar{M}^{\bar{n}}$ and $M^{n}$ at $y_0(=X(x_0))$ and $x_0$ respectively,
    recall that the second fundamental form can be written in
    these local coordinates in the following form
    \begin{equation*}\tag{2.3}
    \begin{split}
h_{ij}^{\alpha}&=\frac{\partial^{2}y^{\alpha}}{\partial
    x^{i}\partial
   x^{j}}-\Gamma^{k}_{ij}\frac{\partial y^{\alpha}}{\partial
   x^{k}}+\bar{\Gamma}^{\alpha}_{\beta\gamma}\frac{\partial y^{\beta}}{\partial x^{i}}
   \frac{\partial y^{\gamma}}{\partial x^{j}}\\
   &=\nabla_{i}\nabla_{j}(y^{\alpha})+\bar{\Gamma}^{\alpha}_{\beta\gamma}\frac{\partial y^{\beta}}{\partial x^{i}}
   \frac{\partial y^{\gamma}}{\partial x^{j}}, \ \ \ \ \ \ \ \ \ \  \mbox{for}\
   \ \
   \alpha=1,2,\cdots, \bar{n},
    \end{split}
    \end{equation*}
where $\nabla_{i}\nabla_{j}(y^{\alpha})$ is the Hessian of
$y^{\alpha}$, which is viewed as a function of $M^{n}$ near $x_0$.
In the following argument, we denote by $\bar{C}_1$ various
constants depending only on $\bar{C}$, $C$ and $\bar{\delta}$.

 Define $f(x)=\bar{d}^2(y_0, X(x))$ on $M^n\cap
X^{-1}(\bar{B}(y_0,\bar{C}_1))$ for some
$\bar{C}_1\leq\bar{\delta}$, then $\nabla_jf=\frac{\partial
f}{\partial y^{\alpha}}\frac{\partial y^{\alpha}}{\partial x^j}$
and the Hessian of $f$ with respect to the metric $g$ on $M^n\cap
X^{-1}(\bar{B}(y_0,\bar{C}_1))$ can be computed as follows
\begin{equation*}\tag{2.4}
\begin{split}
   \nabla_i\nabla_j f
 &=\frac{\partial}{\partial x^i}\nabla_j f-\Gamma_{ij}^k\nabla_k f\\[3mm]
  &=(\frac{\partial^2f}{\partial y^{\alpha}\partial y^{\beta}}
    -\bar{\Gamma}^{\gamma}_{\alpha\beta}\frac{\partial f}{\partial y^{\gamma}})
   \frac{\partial y^{\alpha}}{\partial x^{j}}\frac{\partial y^{\beta}}{\partial x^{i}}
   +\frac{\partial f}{\partial y^{\alpha}}(\frac{\partial^2y^{\alpha}}{\partial x^i\partial x^j}
   -\Gamma_{ij}^k
    \frac{\partial y^{\alpha}}{\partial x^k}
   +\bar{\Gamma}^{\alpha}_{\beta\gamma}
   \frac{\partial y^{\beta}}{\partial x^{i}}\frac{\partial y^{\gamma}}{\partial x^{j}})\\[3mm]
 &=\bar{\nabla}_{\alpha}\bar{\nabla}_{\beta}\bar{d}^2
   \frac{\partial y^{\alpha}}{\partial x^j}\frac{\partial y^{\beta}}{\partial x^i}
   +2\bar{d}\bar{\nabla}_{\alpha}\bar{d}\cdot h_{ij}^{\alpha}.
  \end{split}
\end{equation*}
Using Hessian comparison theorem on $\bar{M}^{\bar{n}}$ and
choosing $\bar{C}_1$ suitable small so that $\bar{d}$ is suitable
small, we get $$ \nabla_i\nabla_j
f\geq\frac{1}{2}g_{ij}\eqno(2.5)$$ on $M^{n}\cap
X^{-1}(\bar{B}(y_0,\bar{C}_1))$.  Now we claim that any closed
geodesic starting and ending at $x_0$ on $(M^{n},g)$ must have
length $\geq 2\bar{C}_{1}.$

We argue by contradiction. Indeed, suppose we have  a closed
geodesic $\gamma:[0,L]\rightarrow M^{n}$ of length $L<2\bar{C}_1,$
$X\circ\gamma$ must be contained in $\bar{B}(y_0,\bar{C}_1)$, then
by (2.5), we have
$$
\frac{d^2}{ds^2}f\circ\gamma(s)=\nabla^{2}f(\dot{\gamma},\dot{\gamma})\geq\frac{1}{2},
\ \ \ \ \ \  \ s\in [0,L].\eqno(2.6)
$$
By the maximum principle, we have
$$
\sup _{s\in[0,L]}f\circ\gamma(s)\leq f\circ\gamma(0),
$$
this implies that $\gamma$ is just a point $\gamma(0).$ The
contradiction  proves the claim.

On the other hand, by the Gauss equation,
$$
R_{ijkl}=\bar{R}_{ijkl}+(h_{ik}^{\alpha}h_{jl}^{\beta}-h_{il}^{\alpha}h_{kj}^{\beta})\bar{g}_{\alpha\beta}(\cdot,0),
$$
we see that
$$
|Rm|\leq \bar{C}+2C^{2}.\eqno(2.7)
$$
Finally, by Klingenberg lemma\cite{CE}, the injectivity radius of
$(M^{n},g)$ at $x_0$ is given by
\begin{eqnarray*}
inj(M^{n},g,x_0)&=& \min\{\text{the conjugate radius at}\  x_0,\\
& &\ \ \ \ \ \ \ \ \ \ \ \ \frac{1}{2}\
\text{the length of the shortest closed geodesic at}\  x_0\}\\
&\geq & \min\{\frac{\pi}{\sqrt{\bar{C}+2C^{2}}},\bar{C}_1 \}.
\end{eqnarray*}
The proof of the theorem is completed.
\end{pf}

 Let $N$ be a
Riemannian manifold,  the distance function $d(y_{1},y_{2})$ can
be regarded as a function on $N\times
  N.$ In the next theorem, we will estimate the
Hessian of the distance function, which is  viewed as the function
of two variables.  The crucial computation of the Hessian was
carried out in \cite{ScY1}.
\begin{thm} Let $N^{n}$ be a complete  Riemannian manifold
  of dimension $n$ satisfying
 $$|{Rm}|\leq{K_0},\ \ \ \ \ \
 inj(N^{n})\geq {i_0}>0.\eqno(2.8)$$
 Let $d(y_{1},y_{2})$ be the distance function
  regarded as a function on $N\times
  N,$ then there is a positive constant $C=C(K_0,i_0)$  such that when $d(y_{1},y_{2})\leq
  \min\{\frac{i_{0}}{2},\frac{1}{4\sqrt{{K_0}}}\}$, we have
  \begin{equation*}\tag{2.9}
  \begin{split}&(i) |\nabla^{2}  d^{2}|(y_1,y_2)\leq C,\\
  &(ii) (\nabla^{2}  d^{2})(X,X)\geq 2|X_1-P_{\gamma}^{-1}X_2|^{2}-C|X|^{2}d^{2}\ \ \ \mbox{for all }
 X\in T_{(y_1,y_2)}N^{n}\times N^{n},
 \end{split}
 \end{equation*}
 where $X=X_1+X_2,$ $X_1\in T_{y_1}N^{n},$ $X_2\in T_{y_2}N^{n},$
 $\nabla$ is the covariant derivative of $N\times N,$
 $\gamma$ is the unique geodesic connecting $y_1$ and $y_2$ in $N^{n}$, and $P_{\gamma}$
  is the parallel translation of $N^{n}$ along $\gamma$.
   \end{thm}
  \begin{pf}
  Set
$\psi(y_{1},y_{2})=d^2_{N^{n}}(y_{1},y_{2}).$ Then $\psi$ is a
smooth function of $(y_1,y_2)$ when $d(y_{1},y_{2})\leq
  \min\{\frac{i_{0}}{2},\frac{1}{4\sqrt{{K_0}}}\}.$
  Now we recall the computation of $Hess(\psi)$ in \cite{ScY1}. For any $(u,v)\in
D=\{(u,v): (u,v)\in N^{n}\times
N^{n},d_{N^{n}}(u,v)\leq\min\{\frac{i_{0}}{2},\frac{1}{4\sqrt{{K_0}}}\}\}\setminus\{(u,u):u\in
N^{n}\} $, let $\gamma_{uv}$ be the minimal geodesic from $u$ to
$v$ and $e_{1}\in T_{u}N^{n}$ be the tangent vector to
$\gamma_{uv}$ at $u$. Then $e_{1}(u,v)$ defines a smooth vector
field on $D$. Let $\{e_{i}\}$ be an orthonormal basis for
$T_{u}N^{n}$ which depends on $u$ smoothly. By parallel
translation of $\{e_{i}\}$ along $\gamma$, we define
$\{\bar{e}_{i}\}$ an orthonormal basis
 for $T_{v}N^{n}$. Thus $\{e_{1},\cdots e_{n},\bar{e}_{1},\cdots \bar{e}_{n}\}$
 is a local frame on $D$. Then for any $X=X_1+X_2\in T_{(u,v)}D$
 with
$$
X_1=\sum_{i=1}^{n}\xi_{i}e_{i} \ \  \mbox{and}\ \
X_2=\sum_{i=1}^{n}\eta_{i}\bar{e}_{i},
$$
  by the formula (16) in
  \cite{ScY1},
\begin{equation*}\tag{2.10}
\begin{split}
\frac{1}{2}Hess(\psi)(X,X)=&\sum_{i=1}^{n}(\xi_{i}-\eta_i)^{2}
+\int_{0}^{r}t\langle\nabla_{e_{1}}V,\nabla_{e_{1}}V\rangle
+\int_{0}^{r}t\langle\nabla_{\bar{e}_{1}}V,\nabla_{\bar{e}_{1}}V\rangle\\
&-\int_{0}^{r}t\langle
R(e_{1},V)V,e_{1}\rangle-\int_{0}^{r}t\langle
R(\bar{e}_{1},V)V,\bar{e}_{1}\rangle,
\end{split}
\end{equation*}
where $V$ is a Jacobi field on geodesic $\sigma$ (connecting
$(v,v)$ to $(u,v)$) and $\bar{\sigma}$ (connecting $(u,u)$ to
$(u,v)$ of length $r=\sqrt{\psi}$) with $X$ as the boundary
values, where $X$ is extended to be a local vector field by
letting its coefficients with respect to $\{e_{1},\cdots
e_{n},\bar{e}_{1},\cdots \bar{e}_{n}\}$ be constant(see
\cite{ScY1}). By the Jacobi equation, we have the estimates
$$|V|\leq C(K_0,i_0)|X|,\ \ \ r|\nabla_{e_{1}}V|\leq
C(K_0,i_0)|X|, \ \ \ \  \ \ \ r|\nabla_{\bar{e}_{1}}V|\leq
C(K_0,i_0)|X|$$  under
  the assumption $d(y_{1},y_{2})\leq
  \min\{\frac{i_{0}}{2},\frac{1}{4\sqrt{{K_0}}}\}$.
  Thus by (2.10)
   we have $$ |Hess(\psi)|\leq C(K_0,i_0),$$ this proves (i).
   Similarly, when $d(y_{1},y_{2})\leq\min\{\frac{i_{0}}{2},\frac{1}{4\sqrt{{K_0}}}\},$ by (2.10), we have
  \begin{eqnarray*}
   \frac{1}{2}Hess(\psi)(X,X)&\geq&\sum_{i=1}^{n}(\xi_{i}-\eta_i)^{2}-\int_{0}^{r}t\langle
R(e_{1},V)V,e_{1}\rangle-\int_{0}^{r}t\langle
R(\bar{e}_{1},V)V,\bar{e}_{1}\rangle
\\
&\geq& \sum_{i=1}^{n}(\xi_{i}-\eta_i)^{2}-C(K_0,i_0)|X|^{2}r^{2}.
  \end{eqnarray*}
This proves (ii). The Theorem is proved.

\end{pf}

 For future applications, in the next
part of this section, we will calculate the equations of
derivatives of general harmonic map flow. Since the MCF is a kind
of harmonic map flow with varying base metrics evolved by MCF, the
formulas computed here are very useful in deriving the higher
derivatives estimates in section 3 and 4. The formulas are of
interest in their own rights. First we fix some notations.

Let $F$ be a map from a Riemannian manifold $(M,g_{ij})$ to
another Riemannian manifold $(N,\bar{g}_{\alpha\beta})$, let
$F^{-1}TN$ be the pull back of the tangent bundle of $N$, we equip
the bundle $(T^{*}M)^{\otimes p}\otimes F^{-1}TN$ the connection
and metric induced from the connections and metrics of $M$ and
$N$. Let $u$ be a section of $(T^{*}M)^{\otimes (p-1)}\otimes
F^{-1}TN$. In local coordinates $\{x^{i}\}$ and $\{y^{\alpha}\}$
of $M$ and $N$ with $y=F(x)$, we have
$|u|^{2}=u^{\alpha}_{i_1i_2\cdots i_{p-1}}u^{\beta}_{j_1j_2\cdots
j_{p-1}}g^{i_1j_1}\cdots g^{i_{p-1}j_{p-1}}\bar{g}_{\alpha\beta}$.
The coefficients of the covariant derivative  $\nabla u$ can be
computed  by the formula
$$
(\nabla u)^{\alpha}_{i_1i_2\cdots i_{p-1}i_p} =\frac{\partial
u_{i_1i_2\cdots i_{p-1}}^{\alpha}}{\partial
x^{i_p}}-\Gamma^{l}_{i_pi_j}u_{i_1i_2\cdots i_{j-1} li_{j+1}\cdots
i_{p-1}}^{\alpha}+\bar{\Gamma}^{\alpha}_{\beta\gamma}\frac{\partial
F^{\beta}}{\partial x^{i_p}}u^{\gamma}_{i_1i_2\cdots i_{p-1}},
$$where $\Gamma$ and $\bar{\Gamma}$ are connection coefficients of
$M$ and $N$ respectively. We can define the Laplacian  of $u$  by
$\triangle u=tr_{g}\nabla^{2}u=g^{ij}(\nabla^{2}u)_{\cdots ij}$.
Recall the Ricci identity $$ (\nabla^{2}u)_{\cdots
ij}^{\alpha}-(\nabla^{2}u)_{\cdots
ji}^{\alpha}=-R_{iji_{m}l}u^{\alpha}_{\cdots
i_{m-1}ki_{m+1}\cdots}g^{kl}+\bar{R}_{\beta\gamma\delta\zeta}\frac{\partial
F^{\beta}}{\partial x^{j}}\frac{\partial F^{\gamma}}{\partial
x^{i}}\bar{g}^{\alpha\delta}u^{\zeta}_{\cdots}. \eqno{(2.11)}
$$ Note that the derivative $\nabla F$ ($\nabla_i
F^{\alpha}=\frac{\partial F^{\alpha}}{\partial x^{i}}$) is a
section of the bundle $T^{*}M\otimes F^{-1}TN$, the higher
derivative $\nabla^{p} F$ is a section of $(T^{*}M)^{\otimes
p}\otimes F^{-1}TN$.

If we have a family of metrics $g_{ij}(\cdot,t)$ on $M$ and a
family of maps $F(\cdot,t)$ from $M$ to $N$, then for each time
$t$, we can still define the bundle $(T^{*}M)^{\otimes p}\otimes
F^{-1}TN$ and define the covariant derivative $\nabla.$  It is a
useful observation that the natural time derivative
$\frac{\partial}{\partial t}$ is not covariant with the metrics.
We define a covariant time derivative $D_t$ as follows. For any
section $u^{\alpha}_{i_1\cdots i_p}$ of $(T^{*}M)^{\otimes
p}\otimes F^{-1}TN$, we define
$$ D_tu^{\alpha}_{i_1\cdots i_p}= \frac{\partial}{\partial
t}u^{\alpha}_{i_1\cdots
i_p}+\bar{\Gamma}^{\alpha}_{\beta\gamma}\frac{\partial
F^{\beta}}{\partial t}u^{\gamma}_{i_1\cdots i_p}.
$$
It is a routine computation which shows that the operator $D_t$ is
covariant.

\begin{prop} Let $M$ be a manifold with a family of
 metrics $g_{ij}(x,t)$, $(N,\bar{g})$ a Riemannian manifold.
 Let  $F(\cdot,t)$ be a solution to the harmonic map flow with respect to the evolving metrics
 $g_{t}$ and $\bar{g}$
   $$
\frac{\partial}{\partial t}F(x,t)=\triangle F(x,t), \ \ \ \ \ \ \
\ \ \text{for}\ x\in M^{n} \text{ and }\  t\geq 0, \eqno{(2.12)}
$$ where
$\triangle F(x,t)$ is the harmonic map Laplacian of $F$ defined by
metrics $g_{ij}(x,t)$ and $\bar{g}$. Then we have

\begin{equation*} \tag{2.13}     \begin{split}
        ( D_t-
     {\triangle}) {\nabla^{k}}{F}&=\sum_{l=0}^{k-1}\nabla^{l}[(R_{M}\ast
g^{-2}+\bar{R}_{N}\ast(\nabla F)^{2}\ast
g^{-1}\ast\bar{g}^{-1})]\ast\nabla^{k-l}F\\&  \ \
+\sum_{l=1}^{k-1}g^{-1}\ast\nabla^{l} \frac{\partial g}{\partial
t}\ast\nabla^{k-l}F,
  \end{split}
            \end{equation*}
            where ${\nabla^{l}}(A\ast B)$ represents
            the linear combinations of ${\nabla^{l}}A\ast
            B$,${\nabla^{l-1}}A\ast
            {\nabla}B$, $\cdots$, $A\ast
            {\nabla^{l}}B$ with universal coefficients.
\end{prop}
\begin{pf} For $k=1$, by direct computation and Ricci identity, we have
\begin{eqnarray*}
\frac{\partial}{\partial
t}{\nabla_{i}}{F^{\alpha}}+\bar{\Gamma}^{\alpha}_{\beta\gamma}F^{\beta}_{i}(\triangle
F)^{\gamma}&=&\nabla_{i}\triangle
F^{\alpha}\\&=&{\triangle}{\nabla_{i}}{F^{\alpha}}
-{R^{l}_{i}}{\nabla_{l}}{F^{\alpha}}+{\bar{R}}^{\alpha}_{\beta\delta\gamma}
{\nabla_{i}}{F^{\beta}}{\nabla_{k}}{F^{\delta}}{\nabla_{l}}{F^{\gamma}}
g^{kl}.\end{eqnarray*}  For $k\geq 2$, we prove by induction.
Since
\begin{eqnarray*}
\frac{\partial}{\partial t}(\nabla^{k}F)^{\alpha}_{i_1\cdots
i_k}&=&\frac{\partial}{\partial x^{i_k}}\frac{\partial}{\partial
t}(\nabla^{k-1}F)^{\alpha}_{i_1\cdots i_{k-1}}-\Gamma^{p}_{i_k
i_l}\frac{\partial}{\partial t}(\nabla^{k-1}F)^{\alpha}_{i_1\cdots
p\cdots i_{k-1}}\\
&
&+\bar{\Gamma}^{\alpha}_{\beta\gamma}F^{\beta}_{i_k}\frac{\partial}{\partial
t}(\nabla^{k-1}F)^{\gamma}_{i_1\cdots i_{k-1}} + (g^{-1}\ast
\nabla \frac{\partial g}{\partial t}\ast
\nabla^{k-1}F)^{\alpha}_{i_1\cdots i_{k}}\\
& &+\frac{\partial}{\partial
y^{\delta}}\bar{\Gamma}^{\alpha}_{\beta\gamma}(\triangle
F)^{\delta}F^{\beta}_{i_k}(\nabla^{k-1}F)^{\gamma}_{i_{1}\cdots
i_{k-1}}+\bar{\Gamma}^{\alpha}_{\beta\gamma}\frac{\partial}{\partial
t}F^{\beta}_{i_k}(\nabla^{k-1}F)^{\gamma}_{i_{1}\cdots i_{k-1}},
\end{eqnarray*}
we have
\begin{eqnarray*}
D_t(\nabla^{k}F)^{\alpha}_{i_1\cdots
i_k}&=&\frac{\partial}{\partial
x^{i_k}}D_t(\nabla^{k-1}F)^{\alpha}_{i_1\cdots
i_{k-1}}-\Gamma^{p}_{i_k
i_l}D_t(\nabla^{k-1}F)^{\alpha}_{i_1\cdots
p\cdots i_{k-1}}\\
&
&+\bar{\Gamma}^{\alpha}_{\beta\gamma}F^{\beta}_{i_k}D_t(\nabla^{k-1}F)^{\gamma}_{i_1\cdots
i_{k-1}} + (g^{-1}\ast \nabla \frac{\partial g}{\partial t}\ast
\nabla^{k-1}F)^{\alpha}_{i_1\cdots i_{k}}\\
& &+\frac{\partial}{\partial
y^{\delta}}\bar{\Gamma}^{\alpha}_{\beta\gamma}(\triangle
F)^{\delta}F^{\beta}_{i_k}(\nabla^{k-1}F)^{\gamma}_{i_{1}\cdots
i_{k-1}}+\bar{\Gamma}^{\alpha}_{\beta\gamma}\frac{\partial}{\partial
t}F^{\beta}_{i_k}(\nabla^{k-1}F)^{\gamma}_{i_{1}\cdots i_{k-1}}\\
&&- \frac{\partial}{\partial
x^{i_k}}[\bar{\Gamma}^{\alpha}_{\beta\gamma}\frac{\partial
F^{\beta}}{\partial t}(\nabla^{k-1}F)^{\gamma}_{i_1\cdots
i_{k-1}}]+\Gamma^{p}_{i_k i_l}
\bar{\Gamma}^{\alpha}_{\beta\gamma}\frac{\partial
F^{\beta}}{\partial t}(\nabla^{k-1}F)^{\gamma}_{i_1\cdots p\cdots
i_{k-1}}\\
& &-\bar{\Gamma}^{\alpha}_{\beta\gamma}
\bar{\Gamma}^{\gamma}_{\delta\xi}F^{\beta}_{i_k}\frac{\partial
F^{\delta}}{\partial t}(\nabla^{k-1}F)^{\xi}_{i_1\cdots
i_{k-1}}+\bar{\Gamma}^{\alpha}_{\beta\gamma}\frac{\partial
F^{\beta}}{\partial t}(\nabla^{k}F)^{\gamma}_{i_1\cdots i_k}.
\end{eqnarray*}
Since
\begin{eqnarray*}
\frac{\partial}{\partial
x^{i_k}}[\bar{\Gamma}^{\alpha}_{\beta\gamma}\frac{\partial
F^{\beta}}{\partial t}(\nabla^{k-1}F)^{\gamma}_{i_1\cdots
i_{k-1}}]&=&\frac{\partial }{\partial
y^{\beta}}\bar{\Gamma}^{\alpha}_{\delta\gamma}F^{\beta}_{i_k}\frac{\partial
F^{\delta}}{\partial t}(\nabla^{k-1}F)^{\gamma}_{i_1\cdots
i_{k-1}}+\bar{\Gamma}^{\alpha}_{\beta\gamma}\frac{\partial}{\partial
x^{i_k}}\frac{\partial F^{\beta}}{\partial
t}(\nabla^{k-1}F)^{\gamma}_{i_1\cdots i_{k-1}}\\
&& +\bar{\Gamma}^{\alpha}_{\beta\gamma}\frac{\partial
F^{\beta}}{\partial t}(\nabla^{k}F)^{\gamma}_{i_1\cdots
i_{k}}+\Gamma^{p}_{i_k i_l}
\bar{\Gamma}^{\alpha}_{\beta\gamma}\frac{\partial
F^{\beta}}{\partial t}(\nabla^{k-1}F)^{\gamma}_{i_1\cdots p\cdots
i_{k-1}}\\&&-\bar{\Gamma}^{\alpha}_{\delta\gamma}
\bar{\Gamma}^{\gamma}_{\beta\xi}F^{\beta}_{i_k}\frac{\partial
F^{\delta}}{\partial t}(\nabla^{k-1}F)^{\xi}_{i_1\cdots i_{k-1}},
\end{eqnarray*}
we have
\begin{eqnarray*}
D_t(\nabla^{k}F)^{\alpha}_{i_1\cdots i_k}&=&[\nabla
D_t(\nabla^{k-1}F)]^{\alpha}_{i_1\cdots i_{k}} + (g^{-1}\ast
\nabla \frac{\partial g}{\partial t}\ast
\nabla^{k-1}F)^{\alpha}_{i_1\cdots i_{k}}\\
& &+\bar{R}^{\alpha}_{\delta\beta\gamma}(\triangle
F)^{\delta}F^{\beta}_{i_k}(\nabla^{k-1}F)^{\gamma}_{i_{1}\cdots
i_{k-1}}.
\end{eqnarray*}
 Combining with Ricci identity
\begin{eqnarray*}
\nabla\triangle\nabla^{k-1}F&=&\triangle
\nabla^{k}F+\nabla[(R_{M}\ast g^{-2}+\bar{R}_{N}\ast(\nabla
F)^{2}\ast
g^{-1}\ast\bar{g}^{-1})\ast\nabla^{k-1}F]\end{eqnarray*}and
induction on $k$, we have
\begin{eqnarray*}
(D_{t}-\triangle)(\nabla^{k}F)&=&g^{-1}\ast\nabla \frac{\partial
g}{\partial t}\ast\nabla^{k-1}F+\bar{R}_{N}\ast\nabla
F\ast\nabla^{2} F\ast\nabla^{k-1}F\ast g^{-1}\ast\bar{g}^{-1}\\&
&+\nabla[(D_{ t}-\triangle)\nabla^{k-1}F]+\nabla[(R_{M}\ast
g^{-2}+\bar{R}_{N}\ast(\nabla
F)^{2}\ast g^{-1}\ast\bar{g}^{-1})\ast\nabla^{k-1}F]\\
&=&\nabla[(D_{t}-\triangle)\nabla^{k-1}F]+\nabla[(R_{M}\ast
g^{-2}+\bar{R}_{N}\ast(\nabla F)^{2}\ast
g^{-1}\ast\bar{g}^{-1})\ast\nabla^{k-1}F]\\& &+g^{-1}\ast\nabla
\frac{\partial
g}{\partial t}\ast\nabla^{k-1}F\\
&=&\sum_{l=0}^{k-1}\nabla^{l}[(R_{M}\ast
g^{-2}+\bar{R}_{N}\ast(\nabla F)^{2}\ast
g^{-1}\ast\bar{g}^{-1})]\ast\nabla^{k-l}F\\&
&+\sum_{l=1}^{k-1}g^{-1}\ast\nabla^{l} \frac{\partial g}{\partial
t}\ast\nabla^{k-l}F.
\end{eqnarray*}
We finish the proof of the proposition.
\end{pf}

\begin{cor}
 Let  $F(\cdot,t)$ be assumed as in proposition 2.3.
    Then we have

\begin{equation*} \tag{2.14}     \begin{split}
        ( \frac{\partial}{\partial t}-
     {\triangle}) |{\nabla^{k}}{F}|^{2}&\leq -2|{\nabla^{k+1}}{F}|^{2}+\langle\sum_{l=0}^{k-1}\{\nabla^{l}[(R_{M}\ast
g^{-2}+\bar{R}_{N}\ast(\nabla F)^{2}\ast
g^{-1}\ast\bar{g}^{-1})]\\&+g^{-1}\ast\nabla^{l} \frac{\partial
g}{\partial
t}\}\ast\nabla^{k-l}F,\nabla^{k}F\rangle+g^{-(k+1)}\frac{\partial
g}{\partial t}\ast(\nabla^{k}F)^{2}\ast \bar{g}.
  \end{split}
            \end{equation*}
\end{cor}
\begin{pf} Since $|\nabla^{k}F|^{2}=(\nabla^{k}F)^{\alpha}_{i_1\cdots
i_{k}}(\nabla^{k}F)^{\beta}_{j_1 \cdots j_{k}}g^{i_1 j_1}\cdots
g^{i_{k} j_{k}}\bar{g}_{\alpha\beta}$, and
\begin{eqnarray*}
\frac{\partial}{\partial
t}|\nabla^{k}F|^{2}&=&2\frac{\partial}{\partial
t}(\nabla^{k}F)^{\alpha}_{i_1 \cdots
i_{k}}(\nabla^{k}F)^{\beta}_{j_1 \cdots j_{k}}g^{i_1 j_1}\cdots
g^{i_{k} j_{k}}\bar{g}_{\alpha\beta}\\&&+\frac{\partial
\bar{g}_{\alpha\beta}}{\partial y^{\delta}}\frac{\partial
F^{\delta}}{\partial t}(\nabla^{k}F)^{\alpha}_{i_1 \cdots
i_{k}}(\nabla^{k}F)^{\beta}_{j_1 \cdots j_{k}}g^{i_1 j_1}\cdots
g^{i_{k} j_{k}}+g^{-(k+1)}\ast\frac{\partial g}{\partial
t}\ast(\nabla^{k}F)^{2}\ast \bar{g}\\
&=& 2D_{t}(\nabla^{k}F)^{\alpha}_{i_1 \cdots
i_{k}}(\nabla^{k}F)^{\beta}_{j_1 \cdots j_{k}}g^{i_1 j_1}\cdots
g^{i_{k} j_{k}}\bar{g}_{\alpha\beta}+g^{-(k+1)}\ast\frac{\partial
g}{\partial t}\ast(\nabla^{k}F)^{2}\ast \bar{g},
\end{eqnarray*}
then (2.14) follows from Proposition 2.3.

\end{pf}

\section{Higher derivative estimates for the mean curvature flow}
\qquad Now we come back to MCF, suppose $X(\cdot,t)$ is a solution
to MCF equation (1.2), $g(\cdot,t)$ is the family of the induced
metrics on $M^{n}$ from $(\bar{M}^{\bar{n}}, \bar{g})$ by
$X(\cdot,t)$, then
$$
\frac{\partial}{\partial
t}g_{ij}=-2H^{\alpha}h^{\beta}_{ij}\bar{g}_{\alpha\beta}.
\eqno{(3.1)}
$$
Note that $\frac{\partial g}{\partial
t}=(\nabla^{2}X)^{2}\ast\bar{g}\ast g^{-1}$ and
$R_{M}=\bar{R}_{\bar{M}}\ast (\nabla
X)^{4}+(\nabla^{2}X)^{2}\ast\bar{g}$. Combining with corollary
2.4, we have
\begin{prop}Let $(\bar{M}^{\bar{n}},\bar{g})$ be a  Riemannian manifold
  of dimension $\bar{n}$.  Let $X_0:M^{n}\rightarrow
 \bar{M}^{\bar{n}}$ be an isometrically immersed  manifold  in $\bar{M}^{\bar{n}}$. Suppose $X(x,t)$
 is a solution of MCF
 on $M^n\times[0,T]$ with  $X_0$ as initial data. Then
\begin{equation*} \tag{3.2}     \begin{split}
        ( \frac{\partial}{\partial t}-
     {\triangle}) |{\nabla^{k}}{X}|^{2}&\leq -2|{\nabla^{k+1}}{X}|^{2}+\langle\sum_{l=0}^{k-1}\nabla^{l}
     [(\nabla^{2}{X})^{2}\ast\bar{g}\ast
g^{-2}+\bar{R}_{\bar{M}}\ast(\nabla X)^{4}\ast g^{-2}\\&\ \ \
\ast\bar{g}\ast\bar{g}^{-1}]
\ast\nabla^{k-l}X,\nabla^{k}X\rangle+g^{-(k+2)}\ast
\bar{g}^{2}\ast(\nabla^{2}X)^{2}\ast(\nabla^{k}X)^{2}.
  \end{split}
            \end{equation*}
\end{prop}

 Now we are ready to derive the higher derivatives
estimates of the second fundamental form of MCF provided that we
have bounded the second fundamental form. Before the deriving of
the higher derivatives estimates, we need to construct a family of
cut-off functions $\xi_k$, which are used also in the next
section. For each integer $k>0$, let $\xi_k$ be a smooth
non-increasing function from $(-\infty,+\infty)$ to $[0,1]$ so
that $\xi_k(s)=1$ for
   $s\in (-\infty,\frac{1}{2}+\frac{1}{2^{k+1}}]$, and $\xi_k(s)=0$ for $s\in[\frac{1}{2}+\frac{1}{2^k},+\infty)$; moreover for
   any $\epsilon>0$ there exists a universal $C_{k,\epsilon}>0$ such
   that
   $$
|\xi^{\prime}_k(s)|+|\xi^{\prime\prime}_k(s)|\leq
C_{k,\epsilon}{\xi_k(s)}^{1-\epsilon}. \eqno{(3.3)}
$$

 \begin{thm}(local estimates) Let $(\bar{M}^{\bar{n}},\bar{g})$ be a complete  Riemannian manifold
  of dimension $\bar{n}$.  Let $X_0:M^{n}\rightarrow
 \bar{M}^{\bar{n}}$ be an isometrically immersed complete manifold  in $\bar{M}^{\bar{n}}$. Suppose $X(x,t)$
 is a solution to the mean curvature flow (1.1)
 on $M^n\times[0,T]$ with  $X_0$ as initial data and with
 bounded second fundamental forms $|h_{ij}^{\alpha}|\leq \bar{C}$ on $[0,T]$.
 Then for any fixed $x_0\in M^{n}$ and any geodesic ball  $B_{0}(x_0,a)$ of radius
$a>0$ of initial metric $g_{ij}$, for any $k\geq 3$, we have
$$
|\nabla^{k}X|(x,t)\leq \frac{C_{k}}{t^{\frac{k-2}{2}}},  \ \ \ \
\mbox{for all}\ \ (x,t)\in B_{0}(x_0,\frac{a}{2})\times[0,T],
\eqno{(3.4)}
$$
where the constant $C_k$ depends on $\bar{C},$ $T,$ $\bar{n},$ $a$
and the bounds of the curvature and its covariant derivatives up
to order $k-1$ of the ambient manifold $\bar{M}$ on its geodesic
ball $B_{\bar{M}}(X_0(x_0),a+1+\sqrt{n}\bar{C}T).$
   \end{thm}
\begin{pf} Since $|\frac{\partial}{\partial t}X|=|H|\leq
\sqrt{n}\bar{C}$, it is not hard to see that under the evolution
of MCF, at any time $t\in [0,T]$, $X_t(B_{0}(x_0,a))$ is contained
in $B_{\bar{M}}(X_{0}(x_0),a+1+\sqrt{n}\bar{C}T).$ For any fixed
$a>0$, $k>0$, we denote by $C_k$ various constants depending only
on $a$, $\bar{C}$, $T$, $\bar{n}$ and the bounds of the curvature
and its covariant derivatives up to order $k-1$ of the ambient
manifold $\bar{M}$ on its ball
$B_{\bar{M}}(X_{0}(x_0),a+1+\sqrt{n}\bar{C}T).$

By Proposition 3.1, we have    \begin{equation*} \tag{3.5}
\begin{split}   ( \frac{\partial}{\partial t}-
     {\triangle}) |{\nabla^{2}}{X}|^{2}& \leq -2|{\nabla^{3}}{X}|^{2}
     +C_2 +C_2|{\nabla^{3}}{X}| \\ &\leq -|{\nabla^{3}}{X}|^{2}
     +C_2
\end{split}
            \end{equation*}
and
\begin{equation*} \tag{3.6}     \begin{split}
        ( \frac{\partial}{\partial t}-
     {\triangle}) |{\nabla^{3}}{X}|^{2}&\leq -2|{\nabla^{4}}{X}|^{2}+C_3(
     |{\nabla^{3}}{X}|^{3}+|{\nabla^{3}}{X}|^{2}+|{\nabla^{3}}{X}|+|\nabla^4X||\nabla^3X|)\\
     &\leq -|{\nabla^{4}}{X}|^{2}+C_3
     |{\nabla^{3}}{X}|^{3}+C_3.
  \end{split}
            \end{equation*}
            Combining (3.5) and (3.6), for any constant
            $A>0$ we have
\begin{equation*} \tag{3.7}     \begin{split}
        ( \frac{\partial}{\partial t}-
     {\triangle}) ((A+|{\nabla^{2}}{X}|^{2})|{\nabla^{3}}{X}|^{2})\leq&  (-|{\nabla^{3}}{X}|^{2}
     +C_3)|{\nabla^{3}}{X}|^{2}+8|{\nabla^{3}}{X}|^{2}|{\nabla^{4}}{X}||{\nabla^{2}}{X}|\\&
     +[-|{\nabla^{4}}{X}|^{2}+C_3
     |{\nabla^{3}}{X}|^{3}+C_3](A+|{\nabla^{2}}{X}|^{2}).
  \end{split}
            \end{equation*}
Since $|{\nabla^{2}}{X}|^{2}$ is bounded by assumption, by
choosing $A$ suitable large, let
$u=(A+|{\nabla^{2}}{X}|^{2})|{\nabla^{3}}{X}|^{2}$ and $v=tu$, we
have
$$( \frac{\partial}{\partial t}-
     {\triangle}) u\leq  -\frac{1}{C_{3}}u^{2}+C_3$$
     and
$$( \frac{\partial}{\partial t}-
     {\triangle}) v\leq  \frac{1}{t}(-\frac{1}{C_{3}}v^{2}+C_3).\eqno{(3.8)} $$
     Now we need a cut-off function technique as in \cite{CZ2}.  Let
$\xi(x)=\xi_3(\frac{{d}_0(x,x_0)}{a})$, where $\xi_3$ is the
cut-off function satisfying (3.3) for $k=3$. Then the function
$\xi(x)$ satisfies
\begin{equation*}\tag{3.9}
\left\{
\begin{split}
     &\xi(x)=1,\ \ \ \ \mbox{ for }x\in B_0(x_0,(\frac{1}{2}+\frac{1}{2^{4}})a),  \\
     &\xi(x)=0, \ \ \ \ \mbox{ for } x \in M\backslash B_0(x_0,a),\\
     &|{\nabla}_0\xi|^2\leq C_3\xi,\\
     &({\nabla}^{2}_0\xi)_{ij}\geq-C_3\xi^{\frac{1}{2}}{g}_{ij}(\cdot,0),
  \end{split}
 \right.
\end{equation*}
where we used the Hessian comparison theorem. Since by Gauss
equation,  the curvature of the initial metric is bounded from
below by a constant, which  depends on  $ \bar{C}$ and the
curvature bound on the ball
$B_{\bar{M}}(X_0(x_0),a+1+\sqrt{n}\bar{C}T)$ of the ambient
manifold.  The last formula holds in the sense of support
functions. Define $\phi(x,t)=\xi(x)v(x,t)$. Then we have
$$( \frac{\partial}{\partial t}-
     {\triangle}) \phi\leq  \frac{1}{t}(-\frac{1}{C_{3}}\xi v^{2}
     -t v\triangle\xi -2t\nabla\xi\cdot\nabla v+C_3\xi).\eqno{(3.10)} $$

Suppose $\phi(x,t)$ achieves its maximum value over
$M^{n}\times[0,T]$ at some point $(x_1,t_1)\in
B(x_0,a)\times(0,T]$, i.e.
$$\phi(x_1,t_1)=\max_{M\times[0,T]}\phi(x,t).$$
 Suppose the point $x_1$
does not lie in the cut-locus of $x_0,$ then
$$\frac{\partial\phi}{\partial t}(x_1,t_1)\geq0,\ \ \ \nabla
v(x_1,t_1)=-\frac{\nabla\xi}{\xi}v,\ \ \ \ \
\triangle\phi(x_1,t_1)\leq0. \eqno{(3.11)}$$
 By (3.10) and (3.11), at $(x_1,t_1)$ we have
$$0\leq-\frac{1}{C_{3}}\xi v^{2}
     -t_1 v\triangle\xi +2 t_1 \frac{|\nabla\xi|^{2}}{\xi}v+C_3\xi.\eqno{(3.12)} $$
Note that the second fundamental form is bounded in
$M^{n}\times[0,T]$, the metrics $g_{ij}(\cdot,t)$ are equivalent.
Since $$\frac{\partial}{\partial t}\Gamma^k_{ij}=(g^{-1}\ast
\nabla \frac{\partial g}{\partial
t})^{k}_{ij}=g^{-2}\ast\bar{g}\ast\nabla^{2}X\ast\nabla^{3}X,$$ we
have \begin{equation*}
\begin{split}
|\Gamma^{k}_{ij}(x_1,t_1)-{\Gamma_0}^{k}_{ij}(x_1)|&\leq
C(\bar{n})\bar{C}\int_0^{t_1}|\nabla^{3}X|dt\\
&\leq C(\bar{n})\bar{C}\int_0^{t_1}(\frac{\phi}{\xi t})^{\frac{1}{2}}(x_1,t)dt\\
&\leq
C_3\frac{\phi(x_1,t_1)^{\frac{1}{2}}}{\xi(x_1)^{\frac{1}{2}}},
\end{split}
\end{equation*}
where we used the fact that $\phi$ achieves its maximum at
$(x_1,t_1).$   Thus at $(x_1,t_1)$, we have
\begin{equation*}
\begin{split}
    -\triangle\xi=&-g^{ij}\nabla_i\nabla_j\xi\\
              =&-g^{ij}({\nabla_0}_i{\nabla_0}_j\xi
                +({\Gamma_0}^k_{ij}-\Gamma_{ij}^k)\nabla_{0k}\xi)\\
           \leq&C_3+
           C_3\frac{\phi(x_1,t_1)^{\frac{1}{2}}}{\xi(x_1)^{\frac{1}{2}}}|\nabla\xi|.
\end{split}
\end{equation*}
Substituting into (3.12), multiplying by $\xi(x_1)$ and combining
with (3.9), we have at $(x_1,t_1)$
\begin{eqnarray*}
0&\leq&-\frac{1}{C_{3}}\xi^{2} v^{2}+
     (C_3+
           C_3\phi(x_1,t_1)^{\frac{1}{2}}\frac{|\nabla\xi|}{\xi^{\frac{1}{2}}}) \xi v+2\frac{|\nabla\xi|^{2}}{\xi}\xi v
           +C_3\xi^{2}\\
&\leq& -\frac{1}{C_{3}}\phi^{2}+C_3
\phi^{\frac{3}{2}}+C_3\phi+C_3.
           \end{eqnarray*}
This implies
$$
\phi(x_1,t_1)\leq C_3,
$$
hence we have
$$
|\nabla^{3}X|\leq \frac{C_3}{t^{\frac{1}{2}}}
$$
on $B_{0}(x_0,(\frac{1}{2}+\frac{1}{2^{4}})a)\times[0,T].$ If
$x_1$ lies on the cut locus of $x_0$, then by applying a standard
support function technique as in \cite{ScY}, the same estimate is
still valid.

For higher derivatives, we prove by induction. Fix $x_0\in M^{n}$,
$a>0$,
 suppose
$$|\nabla^{k}X|\leq\frac{C_k}{t^{\frac{k-2}{2}}},\ \ \ \ k=3,...,m-1,\eqno{(3.13)}$$
on $B_{0}(x_0,(\frac{1}{2}+\frac{1}{2^{k+1}})a)\times[0,T].$  Now
we  prove the estimate for $k=m$.

 By induction hypothesis and
Proposition 3.1, we have
\begin{equation*} \tag{3.14}     \begin{split}
        ( \frac{\partial}{\partial t}-
     {\triangle}) |{\nabla^{m}}{X}|^{2}&\leq -2|{\nabla^{m+1}}{X}|^{2}+\langle\sum_{l=0}^{m-1}\nabla^{l}
     [(\nabla^{2}{X})^{2}\ast\bar{g}\ast
g^{-2}+\bar{R}_{\bar{M}}\ast(\nabla X)^{4}\ast g^{-2}\\&\ \ \ \
\ast\bar{g}\ast\bar{g}^{-1}]
\ast\nabla^{m-l}X,\nabla^{m}X\rangle+g^{-(m+2)}\ast
\bar{g}^{2}\ast(\nabla^{2}X)^{2}\ast(\nabla^{m}X)^{2}\\
&\leq
-2|{\nabla^{m+1}}{X}|^{2}+C_m\sum_{l=0}^{m-1}\{\sum_{l_1+l_2=l}|\nabla^{2+l_1}{X}||\nabla^{2+l_2}{X}|\\&
\ \ \ \ +\sum_{l_1+\cdots
+l_4=l}|{\nabla^{l_1+1}}{X}||{\nabla^{l_2+1}}{X}||{\nabla^{l_3+1}}{X}|
|{\nabla^{l_4+1}}{X}|\}|{\nabla^{m-l}}{X}||{\nabla^{m}}{X}|\\
&\leq-2|{\nabla^{m+1}}{X}|^{2}+C_m[|{\nabla^{m+1}}{X}||{\nabla^{m}}{X}|
+(|{\nabla^{3}}{X}|+1)|{\nabla^{m}}{X}|^{2}\\& \ \ \ \
+t^{-\frac{m-2}{2}}|{\nabla^{m}}{X}|]
\\&
\leq-|{\nabla^{m+1}}{X}|^{2}+\frac{C_m}{t^{\frac{1}{2}}}|{\nabla^{m}}{X}|^{2}+C_mt^{-\frac{m-2}{2}}|{\nabla^{m}}{X}|
  \end{split}
            \end{equation*}
            and
   \begin{equation*} \tag{3.15}     \begin{split}
        ( \frac{\partial}{\partial t}-
     {\triangle}) |{\nabla^{m-1}}{X}|^{2}&\leq -|{\nabla^{m}}{X}|^{2}+\frac{C_{m-1}}{t^{\frac{1}{2}}}
     |{\nabla^{m-1}}{X}|^{2}
     +C_{m-1}t^{-\frac{m-3}{2}}|{\nabla^{m-1}}{X}|\\
     &\leq -|{\nabla^{m}}{X}|^{2}+ \frac{C_{m-1}}{t^{m-3+\frac{1}{2}}}
  \end{split}
            \end{equation*}
on $B_{0}(x_0,(\frac{1}{2}+\frac{1}{2^{m}})a)\times[0,T].$

Define
$$\psi(x,t)=(A+t^{m-3}|\nabla^{m-1}X|^2)|\nabla^m X|^2t^{m-2}$$
for $A$ to be determined later. Combining (3.14) and (3.15), we
have for suitable large $A$ as before
\begin{equation*}\tag{3.16}
\begin{split}
        (\frac{\partial}{\partial
        t}- \triangle) \psi&\leq
        \frac{2m-5}{t}\psi+t^{m-3}|\nabla^m
        X|^2t^{m-2}(-|{\nabla^{m}}{X}|^{2}+
        \frac{C_{m-1}}{t^{m-3+\frac{1}{2}}})\\& \ \ \ \ +t^{m-2}(A+t^{m-3}|\nabla^{m-1}X|^2)(-|{\nabla^{m+1}}{X}|^{2}+\frac{C_m}{t^{\frac{1}{2}}}
        |{\nabla^{m}}{X}|^{2}+C_mt^{-\frac{m-2}{2}}|{\nabla^{m}}{X}|)
        \\&\ \ \ \ +8t^{2m-5}|\nabla^{m-1}X||\nabla^{m}X|^2|\nabla^{m+1}X|\\
        &\leq\frac{2m-5}{t}\psi-\frac{1}{2t}[t^{m-2}|\nabla^m
        X|^2]^{2}+\frac{C_m}{t^{\frac{1}{2}}}[t^{m-2}|\nabla^m
        X|^2]+C_m[t^{m-2}|\nabla^m
        X|^2]^{\frac{1}{2}}\\
        &\leq\frac{1}{t}[-\frac{1}{C_m}\psi^{2}+C_m\psi
        +C_m\psi^{\frac{1}{2}}]\\
        &\leq \frac{1}{t}[-\frac{1}{C_m}\psi^{2}+C_m]
\end{split}
\end{equation*}
on $B_{0}(x_0,(\frac{1}{2}+\frac{1}{2^{m}})a)\times[0,T].$ To
apply the cut-off function technique to (3.16) as before, we note
that by the estimate for $k=3$, we know  that
$$
|\Gamma-\Gamma_0|\leq
C(\bar{n})\bar{C}\int_{0}^{T}|\nabla^{3}X|dt\leq
C_3\int_{0}^{T}\frac{1}{\sqrt{t}}dt\leq C_3.
$$
By calculating the equation of
$\xi_{m}(\frac{d_{0}(x_0,\cdot)}{a})\psi$ using (3.16), and
repeating the same procedure of applying maximum principle as
before, we can prove that
$$
\xi_{m}(\frac{d_{0}(x_0,\cdot)}{a})\psi\leq C_m \ \ \ \ \ \
\mbox{on} \ \ \ \ B_{0}(x_0,a)\times[0,T],
$$
which implies $$ |\nabla^{m}X|(x,t)\leq
\frac{C_{m}}{t^{\frac{m-2}{2}}},  \ \ \ \ \mbox{for all}\ \
(x,t)\in B_{0}(x_0,(\frac{1}{2}+\frac{1}{2^{m+1}})a)\times[0,T].
$$
 We complete the induction step and the theorem is
proved.

\end{pf}

 \begin{cor} Let $(\bar{M}^{\bar{n}},\bar{g})$ be a complete  Riemannian
 manifold satisfying
 $$|\bar{\nabla}^{k}\bar{Rm}|(\cdot)\leq
 \bar{C}, \ \ \ \  \ \ \ \mbox{for  }\ \
 k\leq 2.$$ Let $X_0:M^{n}\rightarrow
 \bar{M}^{\bar{n}}$ be an isometrically immersed complete manifold  in $\bar{M}^{\bar{n}}$. Suppose $X(\cdot,t)$
 is a solution of MCF
 on $M^n\times[0,T]$ with  $X_0$ as initial data and with
 bounded second fundamental forms $|h_{ij}^{\alpha}|\leq \bar{C}$ on $[0,T]$.
 Then there is a constant $C_1$ depending only on $\bar{C},$ $\bar{n}$ and $T$ such that
$$
|\nabla R{m}|(x,t)\leq \frac{C_{1}}{t^{\frac{1}{2}}}, \ \ \
\mbox{for all}\ \ (x,t)\in M^n\times[0,T]. \eqno{(3.17)}
$$ Moreover, for any fixed $x_0\in M^{n}$ and any ball $B_{0}(x_0,a)$ of radius
$a>0$ of initial metric $g_{ij}$, and for any $k\geq 2$, there is
a constant $C_k$ depending only on $a$, $\bar{C}$, $\bar{n}$, $T$
and the bounds of the curvature and its derivatives up to order
$k+1$ of the ambient manifold on its geodesic ball
$B_{\bar{M}}(X_0(x_0),a+1+\sqrt{n}\bar{C}T)$, such that
$$
|\nabla^{k}R{m}|(x,t)\leq \frac{C_{k}}{t^{\frac{k}{2}}},  \ \ \ \
\mbox{for all}\ \ (x,t)\in B_{0}(x_0,\frac{a}{2})\times[0,T].
\eqno{(3.18)}
$$
   \end{cor}
\begin{pf} This follows from Gauss equation and Theorem 3.2.

\end{pf}

\section{Harmonic map flow coupled with mean curvature flow}

\qquad Let $X_t$ be the solution of MCF as in Theorem 1.1,
$g_{ij}(x,t)$ the induced Riemannian metrics.  Let
$f:M^{n}\rightarrow N^{m}$ be a map from $M^{n}$ to  a fixed
Riemanian manifold $(N^{m},\hat{g}_{\alpha\beta}).$ Then the
harmonic map flow coupled with MCF is the following evolution
equation of  maps
\begin{equation*}
\left\{
\begin{split}
 \quad \frac{\partial}{\partial t}f(x,t)&=\triangle f(x,t), \ \ \ \mbox{ for }x \in M^n, t>0,  \\
  f(x,0)&= f(x), \  \  \  \  \   \mbox{ for }x \in M^n ,
  \end{split}
 \right.
\end{equation*}
where the Harmonic map Laplacian $\triangle$ is defined  by using
the metric $g_{ij}(x,t)$ and $\hat{g}_{\alpha\beta}(y)$, i.e.
$$
 \triangle
f^{\alpha}(x,t)=g^{ij}(x,t)\nabla_{i}\nabla_{j}f^{\alpha}(x,t),
   \ \ \ \ \ \ \ \ \ \ $$
   and
   $$
   \nabla_{i}\nabla_{j}f^{\alpha}=\frac{\partial^{2}f^{\alpha}}{\partial
    x^{i}\partial
   x^{j}}-\Gamma^{k}_{ij}\frac{\partial f^{\alpha}}{\partial
   x^{k}}+\hat{\Gamma}^{\alpha}_{\beta\gamma}\frac{\partial f^{\beta}}{\partial x^{i}}
   \frac{\partial f^{\gamma}}{\partial x^{j}}. $$
   Here we use $\{x^{i}\}$ and $\{y^{\alpha}\}$ to denote the
   local
   coordinates of $M^{n}$ and $N^{m}$
   respectively, $\Gamma^{k}_{ij}$ and
   $\hat{\Gamma}^{\alpha}_{\beta\gamma}$ the corresponding Christoffel
   symbols of $g_{ij}$ and $\hat{g}_{\alpha\beta}$.

\vskip 0.5cm  Now we fix a metric $\hat{g}=g(\cdot,T)$ on $M^{n}$,
and let $(N^{m},\hat{g})=(M^{n},\hat{g})$. Note that the ambient
manifold $(\bar{M},\bar{g})$ in Theorem 1.1 satisfies the
assumption of Corollary 3.3. By Corollary 3.3 and Theorem 2.1, we
know that there are positive constants $\hat{C}_1,$ $\hat{\delta}$
depending only on $\bar{C}$, $T$, $\bar{n}$ and $\bar{\delta}$
such that
\begin{equation*} \tag{4.1}
\begin{split}
&|\hat{R}_{N}|+|\hat{\nabla} \hat{R}_{N}|\leq {\hat{C}_{1}},\\
&inj(N,\hat{g})\geq \hat{\delta}>0.
\end{split}
\end{equation*}
 Moreover, by (3.18) of Corollary 3.3, for any fixed $y_0\in N$,  for
any $k\geq 2$, there is a constant $\hat{C}_k$ depending only on
 $\bar{C}$, $\bar{n}$, $T$ and the bounds of the curvature and its
derivatives up to order $k+1$ of the ambient manifold on its ball
$B_{\bar{M}}(X_0(y_0),2e^{\sqrt{n}\bar{C}^{2}T}+1+\sqrt{n}\bar{C}T)$,
such that
$$
|\hat{\nabla}^{k}\hat{R}_{N}|(y)\leq \hat{C_{k}},  \ \ \ \
\mbox{for all}\ \ y\in \hat{B}(y_0,1). \eqno{(4.2)}
$$

In this section, we will establish the existence theorem for the
above harmonic map flow coupled with MCF. More precisely, we will
prove

\begin{thm}
There exists $0<T_0<T$, depending only on $\bar{C}, T, \bar{n},
\bar{\delta}$, such that the harmonic map flow coupled with mean
curvature flow
\begin{equation*}\tag{4.3}
\left\{
\begin{split}
  \frac{\partial}{\partial t}F(x,t)&=\Delta F(x,t), \ \ \ x \in M^n, t>0,  \\
   F(\cdot,0)&=Identity, \  \  \  \  \   x \in M^n
  \end{split}
 \right.
\end{equation*}
has a solution on $M^n\times[0,T_0]$ such that the follwing
estimates hold.  There is a constant $C_2$ depending only on
$\bar{C}$, $\bar{\delta},$ $\bar{n}$ and $T$ such that
$$|\nabla F|+|\nabla^{2} F|\leq C_2. \eqno{(4.4)}$$
For any $k\geq 3$, $B_0(x_1,1)\subset M^n$, there is a constant
$C_k$ depending only on $\bar{C}$,$\bar{\delta},$ $T,$ $\bar{n}$
and $x_1$ such that
$$|\nabla^{k} F|\leq C_{k} t^{-\frac{k-2}{2}}, \ \ \ \mbox{on}\ \
\ B_0(x_1,1)\times[0,T_0].\eqno{(4.5)}
$$
\end{thm}
\vskip0.3cm

We will adapt the strategy of \cite{CZ2} by solving the
corresponding initial-boundary value problem on a sequence of
exhausted bounded domains $D_1\subseteq D_2 \subseteq \cdots $
with smooth boundaries and $D_j\supseteq B_0(x_0,j+1)$,
\begin{equation*}\tag{4.6}
\left\{
\begin{split}
  \frac{\partial}{\partial t}F^{j}(x,t)&=\Delta F^{j}(x,t)  \\
  F^{j}(x,0)&=x \  \  \  \  \   \mbox{for all }x \in D_{j} ,\\
  F^{j}(x,t)&=x \ \ \ \ \ \ \ \ \mbox{for all }x \in \partial D_{j},
\end{split}
\right.
\end{equation*}
and taking a convergent subsequence of $F^{j}$ as
$j\rightarrow\infty$, where $x_0$ is a fixed point in $M^n$.

 First
we need the zero order estimate for the Dirichlet problem (4.6).

\begin{lem}
There exist positive constants $T_1>0$ and $C>0$ such that for any
$j$, if $F^j$ solves problem(4.6) on $\bar{D_{j}}\times[0,T'_1]$
with $T'_1\leq T_1$, then we have
$$\hat{d}(x,F^j(x,t))\leq C\sqrt{t}$$
for any $(x,t)\in D_{j}\times[0,T'_1]$, where $\hat{d}$ is the
distance with respect to the metric $\hat{g}$.
\end{lem}
\begin{pf} \ For simplicity, we drop the superscript $j$. In the
following argument, we denote by $C$ various positive constants
depending only on the constants $\bar{C}$, $\bar{\delta}$, $T$,
and $\bar{n}$ in Theorem 1.1. Note that $\hat{d}(y_{1},y_{2})$ is
the distance function on the target $(M^n,\hat{g})$, which can be
regarded as a function on $M^n \times M^n$ with the product
metric. Let
$\varphi(y_{1},y_{2})=\frac{1}{2}\hat{d}^{2}(y_{1},y_{2})$ and
$\rho(x,t)=\varphi(x,F(x,t))$. We compute
\begin{equation*}
\begin{split}
  (\frac{\partial}{\partial t}-\Delta)\rho
  =&\hat{d}(x,F(x,t))(-\frac{\partial \hat{d}}{\partial y_1^{\alpha}}\Delta Id^{\alpha})
   -g^{ij}\{\frac{\partial^2{\varphi}}{\partial {y_1}^{\alpha}\partial y_1^{\beta}}
    -(\hat{\Gamma}^{\gamma}_{\alpha\beta}\circ Id)
    \frac{\partial\varphi}{\partial y_1^{\gamma}}\}
    \frac{\partial Id^{\alpha}}{\partial x^i}\frac{\partial Id^{\beta}}{\partial x^j}\\
   &-2g^{ij}\frac{\partial^2\varphi}{\partial y_1^{\alpha}\partial y_2^{\beta}}
    \frac{\partial Id^{\alpha}}{\partial x^i}\frac{\partial F^{\beta}}{\partial x^j}
   -g^{ij}\{\frac{\partial^2{\varphi}}{\partial y_2^{\alpha}\partial y_2^{\beta}}
    -(\hat{\Gamma}^{\gamma}_{\alpha\beta}\circ F)
    \frac{\partial\varphi}{\partial y_2^{\gamma}}\}
    \frac{\partial F^{\alpha}}{\partial x^i}\frac{\partial F^{\beta}}{\partial x^j}\\
  =&-\hat{d}\frac{\partial\hat{d}}{\partial y_1^{\alpha}}\triangle Id^{\alpha}-g^{ij}Hess(\varphi)(V_i,V_j),
\end{split}
\end{equation*}
where
$$V_i=\frac{\partial Id^{\alpha}}{\partial x^i}\frac{\partial}{\partial y_1^{\alpha}}
     +\frac{\partial F^{\alpha}}{\partial x^i}\frac{\partial}{\partial y_2^{\alpha}}.$$
 By
Theorem 3.2, there is a constant $C$ depending only on $\bar{C},$
$T$ and $\bar{n}$ such that
$$|\frac{\partial \Gamma}{\partial t}|\leq C|\nabla^3 X|\leq
\frac{C}{\sqrt{t}}.\eqno (4.7)$$ Since
$$\Delta Id=g^{-1}*(\hat{\Gamma}\circ
Id-\Gamma)=g^{-1}*(\Gamma(\cdot,T)-\Gamma(\cdot,t))$$ then  we
have $|\Delta Id|\leq C$ by (4.7),  this implies
$$(\frac{\partial}{\partial t}-\Delta)\rho\leq
C\hat{d}-g^{ij}Hess(\varphi)(V_i,V_j).$$

By (4.1), the curvature of $\hat{g}$ is bounded by some constant
$\hat{K}$, the injectivity radius of $\hat{g}$ has a uniform
positive lower bound $\hat{\delta}$.  We claim that if
$\hat{d}(x,F(x,t))\leq\min\{\hat{\delta}/2,{1}/{4\sqrt{\hat{K}}}\}$,
then
$$g^{ij}Hess(\varphi)(V_i,V_j)\geq\frac{1}{2}|\nabla F|^2-C.$$

Firstly, by Theorem 2.2 (i), we have $|Hess(\varphi)|\leq C$ under
the assumption of the claim. On the other hand, the Hessian
comparison theorem at those points not lying on the cut locus
shows that
\begin{eqnarray*}
  \frac{\partial^2 \varphi}{\partial y_2^{\alpha}\partial y_2^{\beta}}
     -(\hat{\Gamma}^{\gamma}_{\alpha\beta}\circ F)\frac{\partial \varphi}{\partial
     y_{2}^{\gamma}}\geq \frac{\pi}{4}\hat{g}_{\alpha\beta},\\
  \frac{\partial ^2 \varphi}{\partial y_1^{\alpha}\partial  y_1^{\beta}}
     -(\hat{\Gamma}^{\gamma}_{\alpha\beta}\circ Id)\frac{\partial \varphi}{\partial
     y_1^{\gamma}}\geq \frac{\pi}{4}\hat{g}_{\alpha\beta}.
\end{eqnarray*}
Combining the above inequalities, we have
\begin{equation*}
\begin{split}
   g^{ij}Hess(\varphi)(V_i,V_j)
   \geq&\frac{\pi}{4} |\nabla F|^{2}-C|\nabla F|-C\\[3mm]
   \geq&\frac{1}{2}|\nabla F|^2-C,
  \end{split}
\end{equation*}
which proves the claim. Hence when
$\hat{d}(x,F(x,t))\leq\min\{\frac{\hat{\delta}}{2},\frac{1}{4\sqrt{\hat{K}}}\}$,
we have
$$(\frac{\partial}{\partial t}-\Delta)\rho\leq
C\hat{d}-\frac{1}{2}|\nabla F|^2+C.\eqno(4.8)$$ By maximum
principle we have
$$\hat{d}(x,F(x,t))\leq C\sqrt{t}\ \ \ \ \mbox{whenever
}\hat{d}(x,F(x,t))\leq\min\{\frac{\hat{\delta}}{2},\frac{1}{4\sqrt{\hat{K}}}\}.$$
Therefore there exists
$T_1\leq\frac{1}{C^2}\min^2\{\frac{\hat{\delta}}{2},\frac{1}{4\sqrt{\hat{K}}}\}$
such that
$$\hat{d}(x,F(x,t))\leq C\sqrt{t}, \ \ \ \ \mbox{for } t\leq T'_1(\leq T_1),$$
we have proved the lemma.
\end{pf}
 After proving the above lemma, we can apply the standard
parabolic equation theory  to get a local existence for the
initial-boundary value problem (4.6) as follows. This is similar
to \cite{CZ2}, we include the proof here for completeness.
\begin{lem}
There exists a positive constant $T_2 (\leq  T_1)$ depending only
on the dimension $n$, the constants $T_1$ and $C$ obtained in the
previous lemma such that for each $j$, the initial-boundary value
problem (4.6) has a smooth solution $F^{j}$ on
$\bar{D_{j}}\times[0,T_2]$.
\end{lem}
\begin{pf} For an arbitrarily fixed point $\bar{x}$ in $M^{n}$, we
consider the normal coordinates $\{x^i\}$ and $\{y^\alpha \}$ of
the metric $g_{0ij}$ and the metric $\hat{g}_{\alpha \beta}$
respectively around $\bar{x}$. Locally the equation (4.6) is
written as a system of equations
$$ \frac
{\partial {y^{\alpha}}}{\partial t} (x,t) = g^{ij}(x,t)
[\frac{\partial^{2}y^{\alpha}}{\partial
    x^{i}\partial
   x^{j}}-\Gamma^{k}_{ij}(x,t)\frac{\partial y^{\alpha}}{\partial
   x^{k}}
   +\hat{\Gamma}^{\alpha}_{\beta\gamma}(y^1(x,t), \cdots, y^n(x,t))\frac{\partial y^{\beta}}{\partial x^{i}}
   \frac{\partial y^{\gamma}}{\partial x^{j}}].\eqno(4.9)$$
 Note that
$\hat{\Gamma}^{\alpha}_{\beta\gamma}(\bar{x}) = 0$. Since by (4.1)
the curvature of metric $\hat{g}$ and it's first covariant
derivative  are bounded on the whole target manifold, by applying
Corollary 4.12 in \cite{Ha3}, we know that there is some uniform
constant $\hat{C}$ such that if $\hat{d}(y,\bar{x})\leq
\frac{1}{\hat{C}}$, then
$|\hat{\Gamma}^{\alpha}_{\beta\gamma}(y)|\leq
\hat{C}\hat{d}(y,\bar{x}).$ (This fact is  proved essentially in
\cite{Ha3}, although it is not explicitly stated.)  By Lemma 4.2,
$\hat{d}(x,F(x,t))\leq C\sqrt{t}$, we conclude that the
coefficients of the quadratic terms on the RHS of (4.9) can be as
small as we like provided $T_2>0$ sufficiently small (independent
of $\bar{x}$ and $j$).

Now for fixed $j$, we consider the corresponding parabolic system
of the difference of the map $F^j$ and the identity map. Clearly
the coefficients of the quadratic terms of the gradients are also
very small. Thus, whenever (4.9) has a solution on a time interval
$[0,T'_2]$ with $T'_2 \leq T_2$, we can argue exactly as in the
proof of Theorem 6.1 in Chapter VII of the book \cite{LSU} to
bound the norm of $\nabla F^{j}$ on the time interval $[0,T'_2]$
by a positive constant depending only on $g_{0ij}$, and
$\hat{g}_{\alpha\beta}$ over the domain $D_{j+1}$, the
$L^{\infty}$ bound of $F^j$ obtained in the previous lemma, and
the boundary $\partial D_j$. Hence by the same argument as in the
proof of Theorem 7.1 in Chapter VII of the book \cite{LSU}, we
deduce that the initial-boundary value problem (4.9) has a smooth
solution $F^{j}$ on $\bar{D_{j}}\times[0,T_2]$.

\end{pf}

To get a convergent sequence of $F^j$, we need the following
uniform estimates.

\begin{lem}
There exists a positive constant $T_3$, $0<T_3\leq T_2$,
independent of j, such that if $F^j$ solves
\begin{equation*}
\left\{
\begin{split}
   \frac{\partial}{\partial t}F^{j}(x,t)
         &=\Delta F^{j}(x,t)\ \ \ \ \mbox{on }D_{j}\times[0,T_3],\\[3mm]
   F^{j}(x,0)&=x     \ \ \ \ \ \ \ \ \mbox{on }D_j.
\end{split}
\right.
\end{equation*}
Then for any $B_0(x_1,1)\subset D_{j},$ there is a positive
constant $C=C(\bar{C},\bar{\delta},\bar{n},T)$ such that
$$|\nabla F^{j}|+|\nabla^{2} F^{j}|\leq C $$
on $B_{0}(x_1,\frac{1}{2})\times [0,T_3],$ and for any $k\geq3$
there exist constants $C_k=C(k,\bar{C},\bar{\delta},T,\bar{n},
x_1)$ satisfying
$$ |\nabla^{k} F^{j}|\leq C_k t^{-\frac{k-2}{2}}$$
on $B_{0}(x_1,\frac{1}{2})\times [0,T_3].$
\end{lem}
\begin{pf} We drop the superscript $j$. We denote by $C$ various
constants depending only on $\bar{C}$, $\bar{\delta},$ $T$,
$\bar{n}.$ We first estimate $|\nabla F|$. By Corollary 2.4, we
have
\begin{equation*}
\begin{split}
   (\frac{\partial}{\partial t}-{\triangle}) |{\nabla}{F}|^{2}
  &\leq -2|{\nabla}^2{F}|^{2}+\langle([R_{M}\ast g^{-2}
   +\hat{R}_{N}\ast(\nabla F)^{2}\ast g^{-1}\ast\hat{g}^{-1}]\\
  &+g^{-1}\ast\frac{\partial g}{\partial t})\ast\nabla F,\nabla F\rangle
   +g^{-2}\frac{\partial g}{\partial t}\ast(\nabla F)^{2}\ast \hat{g}.
\end{split}
\end{equation*}
Note that $\frac{\partial g}{\partial
t}=(\nabla^2X)^2*\bar{g}*g^{-1},$  $R_M=\bar{R}_{\bar{M}}*(\nabla
X)^4+(\nabla^2X)^2*\bar{g}$, the second fundamental form
$\nabla^2X$ and curvature $\bar{R}_{\bar{M}}$ are bounded by
assumption, we know that $|\frac{\partial g}{\partial t}|$ and
$|R_M|$ are bounded. The above formula gives
$$\frac{\partial}{\partial t}|\nabla F|^2\leq\Delta|\nabla
F|^2-2|\nabla^2 F|^2+C|\nabla F|^2+C|\nabla F|^4.\eqno (4.10)$$ On
the other hand, we know from (4.8) that
$$\frac{\partial}{\partial t}\rho\leq\Delta\rho-\frac{1}{2}|\nabla
F|^2+C,$$ where $\rho(x,t)=\frac{1}{2}\hat{d}^{2}(x,F(x,t)).$ For
any $a>0$ to be determined later, we compute
\begin{equation*}
\begin{split}
   \frac{\partial}{\partial t}[(a+\rho)|\nabla F|^2]
   \leq&\Delta[(a+\rho)|\nabla F|^2]-2\nabla\rho\cdot\nabla|\nabla F|^2\\[3mm]
       &-2(a+\rho)|\nabla^2F|^2+C(a+\rho)|\nabla F|^2+C(a+\rho)|\nabla F|^4\\[3mm]
       &-\frac{1}{2}|\nabla F|^4+C|\nabla F|^2.
\end{split}
\end{equation*}
Since
\begin{equation*}
\begin{split}
   -2\nabla\rho\cdot\nabla|\nabla F|^2\leq
   C\hat{d}(|\nabla F|+|\nabla F|^2)|\nabla^2F|\leq
   C(|\nabla F|^2+|\nabla F|^4)\hat{d}+C\hat{d}|\nabla^2F|^2
\end{split}
\end{equation*}
and $\hat{d}(\cdot,F(\cdot,t))\leq C\sqrt{t}$, by taking
$a=\frac{1}{8C}$ and $T_3$ suitable small, we have
\begin{equation*}
\begin{split}
   \frac{\partial}{\partial t}[(a+\rho)|\nabla F|^2]
   \leq&\triangle[(a+\rho)|\nabla F|^2]
        -\frac{1}{8C}|\nabla^2 F|^2-\frac{1}{4}|\nabla F|^4+C
\end{split}
\end{equation*}
for $t\leq T_3.$ Let $u=(a+\rho)|\nabla F|^2$, then
$$\frac{\partial u}{\partial
t}\leq\Delta u-\frac{1}{C}u^2+C\eqno(4.11)$$ for $t\leq T_3.$ Let
$\xi(x)=\xi_1(d_0(x_1,x))$ be a cut-off function, where $\xi_1$ is
the nonincreasing smooth function in (3.3) supported in $[0,1)$
and equal to $1$ in $[0,\frac{3}{4}]$. Note that at $t=0,$ $u=a
g^{ij}(\cdot,0)g_{ij}(\cdot,T)\leq C$. Then by computing the
equation of $\xi u$ and applying the maximum principle as before,
we have
$$
\xi u(x,t)\leq C  \ \ \ \ \ \ \ \ \mbox{on} \ \ \ \ M^{n}\times
[0,T_3],
$$
this implies
$$
|\nabla F|\leq C \ \ \ \mbox{on} \ \ B_{0}(x_1,\frac{3}{4})\times
[0,T_3].
$$
\vskip0.3cm

We now estimate $|\nabla^2 F|$. By Corollary 2.4 again
  \begin{eqnarray*}
        ( \frac{\partial}{\partial t}-
     {\triangle}) |{\nabla^{2}}{F}|^{2}
     &\leq &-2|{\nabla^{3}}{F}|^{2}+\langle\sum_{l=0}^{1}\{\nabla^{l}[(R_{M}\ast
g^{-2}+\hat{R}_{N}\ast(\nabla F)^{2}\ast
g^{-1}\ast\hat{g}^{-1})]\\&&+g^{-1}\ast\nabla^{l} \frac{\partial
g}{\partial
t}\}\ast\nabla^{2-l}F,\nabla^{2}F\rangle+g^{-(3)}\frac{\partial
g}{\partial t}\ast(\nabla^{2}F)^{2}\ast \hat{g},
\end{eqnarray*}
and by (3.4),(3.17),(4.1), we know $\sqrt{t}|\nabla \frac{\partial
g}{\partial t}|+\sqrt{t}|\nabla R_{M}|+|\hat{\nabla}\hat{R}_{N}|\leq
C$, and
$$\frac{\partial}{\partial t}|\nabla^2 F|^2\leq\Delta|\nabla^2
F|^2-2|\nabla^3
F|^2+C|\nabla^2F|^2+\frac{C}{\sqrt{t}}|\nabla^2F|\eqno(4.12)$$
 on  $ B_{0}(x_1,\frac{3}{4})\times
[0,T_3].$ This implies
$$\frac{\partial}{\partial t}|\nabla^2 F|\leq\Delta|\nabla^2
F|+C|\nabla^2F|+\frac{C}{\sqrt{t}}.\eqno(4.13)$$ By (4.10) we have
$$\frac{\partial}{\partial t}|\nabla F|^2\leq\Delta|\nabla
F|^2-2|\nabla^2F|^2+C.$$ Let $$u=|\nabla^2F|+|\nabla
F|^2-2C\sqrt{t}+2C\sqrt{T},$$ then
$$\frac{\partial}{\partial t}u\leq\Delta u-u^2+C  \ \ \ \ \mbox{ on}\ \ \   B_{0}(x_1,\frac{3}{4})\times
[0,T_3].\eqno(4.14)$$ Define the cutoff function
$\xi(x)=\xi_2({d_0(x_1,x)})$. Note that at $t=0,$
$|\nabla^{2}F|=|\Gamma_0-\hat{\Gamma}|\leq C$, then
$u\mid_{t=0}\leq C$. Using the similar maximum principle argument
as before, we get
$$\xi u\leq C\ \ \ \ \ \ \mbox{on }
B_0(x_1,\frac{1}{2}+\frac{1}{2^2})\times[0,T_3],$$ which implies
$$|\nabla^2 F|\leq C\ \ \ \ \ \mbox{on }
B_0(x_1,\frac{1}{2}+\frac{1}{2^3})\times[0,T_3].$$

To derive the higher derivative estimates we prove by induction on
$k$. We denote by $C_k$ various constants,  depending only on
$\bar{C}$, $T$, $\bar{\delta}$, $\bar{n}$,  and the bounds of the
ambient manifold $\bar{M}$ curvature and its covariant derivatives
up to order $k$  on its ball $B_{\bar{M}}(X_0(x_1),C)$ for
suitable $C$.

Now suppose we have proved
$$|\nabla^l F|\leq\frac{C_l}{t^{\frac{l-2}{2}}},\ \ \ \
l=2,,...,k-1\eqno(4.15)$$ on
$B_{0}(x_1,(\frac{1}{2}+\frac{1}{2^{k}}))\times[0,T_3]$. By
Corollary 2.4, Theorem 3.2, Corollary 3.3  and using (4.15), we
get
$$\frac{\partial}{\partial t}|\nabla^k
F|^2\leq\Delta|\nabla^kF|^2-2|\nabla^{k+1}F|^2+C_k|\nabla^kF|^2
+\frac{C_k}{t^{\frac{k-1}{2}}}|\nabla^kF|,\eqno(4.16)$$ which
implies
$$\frac{\partial}{\partial t}|\nabla^k
F|\leq\Delta|\nabla^kF|+C_k|\nabla^kF|+\frac{C_k}{t^{\frac{k-1}{2}}}.\eqno(4.17)$$
 We also have
$$\frac{\partial}{\partial t}|\nabla^{k-1}
F|^2\leq\Delta|\nabla^{k-1}F|^2-2|\nabla^{k}F|^2
+\frac{C_{k-1}}{t^{k-\frac{5}{2}}}.\eqno(4.18)$$ Let
$$u=t^{\frac{k-2}{2}}|\nabla^k F|+t^{k-3}|\nabla^{k-1}F|^2.$$
By combining (4.17) and (4.18), we obtain
$$\frac{\partial}{\partial t}u\leq\Delta u-\frac{1}{t}(u^2+C_k) \eqno(4.19)$$
 on
$B_{0}(x_1,(\frac{1}{2}+\frac{1}{2^{k}}))\times[0,T_3].$ Using the
cutoff function $\xi(x)=\xi_k(d_0(x_1,x))$, (4.19) and applying
maximum principle as before,  we conclude with
$$|\nabla^k F|\leq\frac{C_k}{t^{\frac{k-2}{2}}}\ \ \ \ \ \ \
\mbox{on }
B_0(x_1,(\frac{1}{2}+\frac{1}{2^{k+1}}))\times[0,T_3].$$ Therefore
we complete the proof of Lemma 4.4.
\end{pf}

 \textbf{Proof of Theorem 4.1}

  Now we combine the above
three lemmas to prove Theorem 4.1. We have known that there is a
$T_3>0$ such that for each $j$, the equation
\begin{equation*}
\left\{
\begin{split}
  \frac{\partial}{\partial t}F^{j}(x,t)
             &=\Delta F^{j}(x,t)  \\
   F^{j}(x,0)&= x \  \  \  \  \   \mbox{for all }x \in D_{j} ,\\
   F^{j}(x,t)&= x \ \ \ \ \ \ \ \ \mbox{for all }x \in \partial D_{j}
\end{split}
\right.
\end{equation*}
has a smooth solution $F^{j} $ on $\bar{D}_j\times[0,T_3]$. Since
$D_{j}\supset B_{0}(x_0,j+1)$, by choosing any $x_1\in
B_{0}(x_0,j)$ in Lemma 4.4 we have
$$|\nabla F^j|+|\nabla^2 F^j|\leq C$$
on $B_0(x_0,j)\times[0,T_3],$ where $C$ depends only on $\bar{C}$,
$\bar{n},$  $\bar{\delta}$, $T$. Moreover for any $x_1\in
B_0(x_0,j)$, $k\geq 3$, there is a $C_k$ depending on $\bar{C}$,
$\bar{\delta}$, $T$, $\bar{n}$ and $x_1$ such that
$$|\nabla^k F^j|(x_1,t)\leq C_k t^{-\frac{k-2}{2}}.$$  Then we can
take a convergent subsequence of $F^j$ (as $j\rightarrow \infty$)
to get the desired $F$ with the desired estimates. So the proof of
Theorem 4.1 is completed.

$\hfill\Box$

For later purpose, now we need to derive some preliminary estimate
of $g_{ij}(x,t)$ with respect to $F^{*}\hat{g}$.  Let
$\hat{g}_{ij}=(F^{*}\hat{g})_{ij}$.

\begin{prop}
Under the assumption of Theorem 4.1, there exist $0<T_4\leq T_3$
and $C>0$ depending only on $\bar{C}$, $\bar{n},$  $\bar{\delta}$
and $T$ such that for all $(x,t)\in M^n\times [0,T_4]$, we have
$$ \frac{1}{C}\hat{g}_{ij}(x,t)\leq g_{ij}(x,t)\leq C\hat{g}_{ij}(x,t).\eqno(4.20)$$
\end{prop}
\begin{pf} Note that $|\nabla F|^2=\hat{g}_{ij}g^{ij}\leq C$, which
implies $\hat{g}_{ij}(x,t)\leq Cg_{ij}(x,t)$. For the reverse
inequality, since the curvature of $g_{ij}(\cdot,t)$ is bounded,
we compute the equation of $\hat{g}_{ij}(x,t)$ on the domain,
\begin{equation*}\tag{4.21}
\begin{split}
 \frac{\partial}{\partial t}\hat{g}_{ij}
     =&\Delta\hat{g}_{ij}-R_{ik}F^{\alpha}_{l}F^{\beta}_{j}\hat{g}_{\alpha\beta}g^{kl}
     -R_{jk}F^{\alpha}_{l}F^{\beta}_{i}\hat{g}_{\alpha\beta}g^{kl}
       +2\hat{R}_{\alpha\beta\gamma\delta}F^{\alpha}_{i}F^{\beta}_{k}
       F^{\gamma}_{j}F^{\delta}_{l}g^{kl}-2\hat{g}_{\alpha\beta}
       F^{\alpha}_{ki}F^{\beta}_{lj}g^{kl}\\[3mm]
  \geq&\Delta\hat{g}_{ij}- R_{ik}\hat{g}_{jl}g^{kl}- R_{jk}\hat{g}_{il}g^{kl}
       -C|\nabla F|^{2}\hat{g}_{ij}-2|\nabla^{2}F|^{2}g_{ij}\\[3mm]
  \geq&\Delta\hat{g}_{ij}-Cg_{ij}.
\end{split}
\end{equation*}
 Note that for suitable large constant $C$, we have
$$
\frac{\partial}{\partial t}g_{ij}\leq Cg_{ij}, \ \ \ \ \ 0<t<T,
$$
and    $ \hat{g}_{ij}\geq\frac{1}{C}g_{ij}$ at time $0.$ Thus  for
$t\leq 1/C^{3},$ we have
\begin{equation*}\tag{4.22}
\begin{split}
 (\frac{\partial}{\partial t}-\triangle)
 (\hat{g}_{ij}+(C^{2}t-\frac{1}{C})g_{ij})
   \geq&[-C+C^{2}+C(C^{2}t-\frac{1}{C})]g_{ij}\geq 0.
   \end{split}
\end{equation*}
 Note that
$$(\hat{g}_{ij}+(C^{2}t-\frac{1}{C})g_{ij})\mid_{t=0}\geq0.$$
  Since  $|\nabla^{2}X|+\sqrt{t}|\nabla^{3}X|\leq C$ and the curvature is
bounded,  then there is a smooth proper function $\varphi$ with
$\varphi(x)\geq 1+d_0(x_0,x),$ $|\nabla
\varphi|+|\nabla^{2}\varphi|\leq C$. So Hamilton's maximum principle
for tensors on complete manifolds is applicable, we get
$$\hat{g}_{ij}+(C^{2}t-\frac{1}{C})g_{ij}\geq0     \ \ \ \ \mbox{for}  \ \ \ \ \   t\leq\min\{T_3,C^{-3}\},$$
which implies
$$g_{ij}\leq2C\hat{g}_{ij}$$
for $t\leq T_4=\min\{T_3,1/2C^{3}\}$.

 The proof of the proposition is completed.\end{pf}

\vskip .2cm As a consequence, we know that the solution of the
harmonic map flow coupled with the MCF is a family of
diffeomorphisms.

\begin{cor}Let $F(x,t)$ be assumed as in the previous  proposition. Then
$F(\cdot,t)$ are diffeomorphisms from $M$ to $N$ for all $t\in
[0,T_4]$.
\end{cor}
\begin{pf} Note that (4.20) implies that $F$ are local diffeomorphisms.
For any $x_1\neq x_2$, we claim that $ F(x_1,t)\neq F(x_2,t)$ for
all $t\in[0,T_ 4]$. Suppose not, then there is the first time
$t_{0} > 0$ such that $ F(x_1,t_{0})= F(x_2,t_{0})$. Choose small
$\sigma>0$ so that there exist a neighborhood $\hat{O}$ of
$F(x_1,t_0)$ and a neighborhood $O$ of $x_1$ such that
$F^{-1}(\cdot,t)$ is a diffeomorphism from $\hat{O}$ to $O$ for
each $t\in[t_0-\sigma,t_0]$, and let $\hat{\gamma}$ be a shortest
geodesic( parametrized by arc length) on the target (with respect
to the metric $\hat{g}$) with $\hat{\gamma}(0)=F(x_1,t)$,
$\hat{\gamma}(l)=F(x_2,t)$ and $\hat{\gamma}\subset\hat{O}$. We
compute
$$
 \frac{\partial}{\partial t}\hat{d}(F(x_1,t),F(x_2,t))=\langle
 V(F(x_2,t)),{\hat{\gamma}}^{\prime}(l)\rangle_{\hat{g}}-\langle
 V(F(x_1,t)),{\hat{\gamma}}^{\prime}(0)\rangle_{\hat{g}}, \eqno{(4.23)}
$$
 where $ V(F(x,t))=\frac{\partial}{\partial
{t}}F(x,t)$. Now we pull back everything by $F^{-1}$ to $O$,
\begin{eqnarray*}
 \frac{\partial}{\partial t}\hat{d}(F(x_1,t),F(x_2,t))&=&\langle
 P_{-\hat{\gamma}}V-V,\hat{\gamma}^{\prime}(0)\rangle_{F^{*}\hat{g}}\\
 &\geq&- \sup_{x\in F^{-1}\hat{\gamma}}|\hat{\nabla}V|(x,t)
 \hat{d}(F(x_1,t),F(x_2,t)),
\end{eqnarray*}
where $P_{\hat{\gamma}}$ is the parallel translation along
$F^{-1}\hat{\gamma}$ using the connection defined by
$F^{*}\hat{g}.$ Since  $$\hat{\nabla}_{k}
V^{l}=\nabla_{k}V^{\alpha}\frac{\partial x^{l}}{\partial
 y^{\alpha}},$$
where $\nabla_{k}V^{\alpha}$ is the covariant derivative of the
section $V^{\alpha}$ of the bundle $F^{-1}TN.$ Thus by (4.20) in
proposition 4.5, we have
$$|\hat{\nabla}_{k}
V^{l}|=[{\nabla}_{k} V^{\alpha}{\nabla}_{l}
V^{\beta}\hat{g}_{\alpha\beta}\hat{g}^{kl}]^{\frac{1}{2}}\leq
C|\nabla^{3}F|\leq \frac{C}{\sqrt{t}},$$
 where the constant $C$ depends on the $x_1$ and $x_2$ and is
 independent of $t$ by (4.5) of Theorem 4.1.
 Therefore, for $t\in [t_0-\sigma,t_0]$, we have
 $$\hat{d}(F(x_1,t),F(x_2,t))\leq e^{C(\sqrt{t_0}-\sqrt{t_0-\sigma})
 }\hat{d}(F(x_1,t_0),F(x_2,t_0))=0,$$ which contradicts with the choice of $t_0$.
 The corollary is proved.
 \end{pf}

%%%%%%%%%%%%%%%%%%%%%%%%%%%%%%%%%%%%%%%%%%%%%%%%%%%%%%%%%%%%%%%%%%%%%%%%%%%%%%%%%%%%%%%%%%
\section{Mean-De Turck flow}
\qquad  From the previous section, we know that the harmonic map
flow coupled with MCF with identity as initial data has a short
time solution $F(x,t)$ which maintains being a diffeomorphism with
good estimates. Let $\bar{X}=X\circ F^{-1}$ be a family of maps
defined from $(N,\hat{g}_{\alpha\beta})$ to $\bar{M}^{\bar{n}}$,
then $\bar{X}$ satisfies the following mean De turck flow
\begin{equation*} \tag{5.1}
\begin{split}
 \frac{\partial}{\partial t}\bar{X}
  = g^{\alpha\beta} \hat{{\nabla}}_{\alpha} \hat{\nabla}_{\beta}\bar{X}\
 \ \ \mbox{ for }y \in N,
  \end{split}
\end{equation*}
where $g^{\alpha\beta}$ is the inverse matrix of
$g_{\alpha\beta}(\cdot,t)=((F^{-1})^{\ast}g(\cdot,t))_{\alpha\beta},$
 $\hat{\nabla}$ is the covariant derivative with respect to
$\hat{g}_{\alpha\beta}$.  We denote the local coordinates of
$\bar{M}$ by $\{\bar{z}^{\bar{\alpha}}\}$. It is not hard to see
$$
g_{\alpha\beta}(y,t)=g_{ij}(x,t)\frac{\partial x^{i}}{\partial
y^{\alpha}}\frac{\partial x^{j}}{\partial
y^{\beta}}=\bar{g}_{\bar{\alpha}\bar{\beta}} \frac{\partial
X^{\bar{\alpha}}}{\partial x^{i}}\frac{\partial
X^{\bar{\beta}}}{\partial x^{j}}\frac{\partial x^{i}}{\partial
y^{\alpha}}\frac{\partial x^{j}}{\partial y^{\beta}} =
\frac{\partial \bar{X}^{\bar{\gamma}}}{\partial
y^{\alpha}}\cdot\frac{\partial \bar{X}^{\bar{\delta}}}{\partial
y^{\beta}}\bar{g}_{\bar{\gamma}\bar{\delta}}(\bar{X}(y,t)),\eqno(5.2)$$
this implies that the metric $g_{\alpha\beta}(y,t)$ is just the
induced metric from the ambient space by the map $\bar{X}$.
Since$$
\hat{\Gamma}^{\gamma}_{\alpha\beta}(y)-\Gamma^{\gamma}_{\alpha\beta}(y,t)=(\nabla^{2}
F)^{\gamma}_{ij}\frac{\partial x^{i}}{\partial
y^{\alpha}}\frac{\partial x^{j}}{\partial y^{\beta}},$$ we have
\begin{equation*}\tag{5.3}
\begin{split}
&\frac{1}{C}\hat{g}_{\alpha\beta}(y)\leq g_{\alpha\beta}(y,t) \leq
C\hat{g}_{\alpha\beta}(y),\\ &
|\hat{\Gamma}^{\gamma}_{\alpha\beta}(y)-\Gamma^{\gamma}_{\alpha\beta}(y,t)|\leq
C,
\end{split}
\end{equation*}
 by Theorem
4.1 and Proposition 4.5.

 Let $X_1$ and $X_2$ be  two solutions of MCF  with bounded second fundamental form and with the same initial
value $X_0$ assumed as in the Theorem 1.1. Let $g_{ij}^{1}(x,t)$ and
$g_{ij}^{2}(x,t)$ be the corresponding induced metrics. As in
section 4, we solve the harmonic map flows coupled with MCF with the
same target $(M^{n},\hat{g}_{\alpha\beta})$ where
$\hat{g}={g}^{1}(T)$ respectively
\begin{equation*}\tag{5.4}
\left\{
\begin{split}
 \frac{\partial}{\partial t}F_1&=\Delta_{g^{1},\hat{g}} F_1  \\
  F_1\mid_{t=0}&= \mbox{Identity} \  \  \  \  \   \mbox{on } M^n ,
  \end{split}
 \right.
\end{equation*}
 and
\begin{equation*}\tag{5.5}
\left\{
\begin{split}
  \frac{\partial}{\partial t}F_2&=\Delta_{g^{2},\hat{g}}  F_2  \\
  F_2\mid_{t=0}&= \mbox{Identity} \  \  \  \  \   \mbox{on }\  M^n ,
  \end{split}
 \right.
\end{equation*}
where $\Delta_{g^{k},\hat{g}} $ is the harmonic map Laplacian
defined by the metric $g_{ij}^{k}(x,t)$ and
$\hat{g}_{\alpha\beta}$ for $k=1,2$ respectively. By section 4, we
obtain two solutions $F_1(x,t)$ and $F_2(x,t)$ such that Theorem
4.1 holds with $F=F_1$ and $F=F_2$. Corollary 4.6 says that
$F_1(x,t)$ and $F_2(x,t)$ are diffeomorphisms for any $t\in
[0,T_4].$ Let
${g_1}_{\alpha\beta}(y,t)=((F_{1}^{-1})^{*}g^{1})_{\alpha\beta}(y,t)$
and
${g_2}_{\alpha\beta}(y,t)=((F_{2}^{-1})^{*}g^{2})_{\alpha\beta}(y,t).$
 Then $\bar{X}_1=X_1\circ F_1^{-1}$ and $\bar{X}_2=X_2\circ
F_2^{-1}$ are two solutions to the mean-De Turck flow (5.1) with
the same initial value $X_0$,
\begin{equation*} \tag{5.6}
\left\{
\begin{split}
 \frac{\partial}{\partial t}\bar{X}_1
  &= g_{1}^{\alpha\beta} \hat{{\nabla}}_{\alpha} \hat{\nabla}_{\beta}\bar{X}_1, \
 \ \ \mbox{on}\  M^n\times[0,T_4],\\
\bar{X}_1\mid_{t=0}&=X_0, \ \ \ \ \ \ \ \ \ \ \ \ \mbox{ on }  M^n,
  \end{split}
  \right.
\end{equation*}
\begin{equation*} \tag{5.7}
\left\{
\begin{split}
 \frac{\partial}{\partial t}\bar{X}_2
  &= g_{2}^{\alpha\beta} \hat{{\nabla}}_{\alpha} \hat{\nabla}_{\beta}\bar{X}_2,\
 \ \ \mbox{ on } M^n\times[0,T_4],\\
\bar{X}_2\mid_{t=0}&=X_0, \ \  \ \ \ \ \ \  \ \ \ \ \ \mbox{ on }
M^n,
  \end{split}
  \right.
\end{equation*}
 where  $g_{1\alpha\beta}$ and
$g_{2\alpha\beta}$  are the corresponding induced metrics from the
target $(\bar{M}^{\bar{n}},\bar{g}_{\bar{\alpha}\bar{\beta}})$ by
the maps $\bar{X}_{1}$ and $\bar{X}_{2}$ by (5.2).

\vskip 0.2cm
\begin{prop} Under the assumptions of Theorem 1.1, there is some
$T_5>0$ depending only on  $\bar{C}$, $\bar{\delta}$, $T$ and
$\bar{n}$ such that
$$\bar{X}_1(y,t)=\bar{X}_2(y,t)\ \ \ \ \ on\ M^n\times[0,T_5]$$
for the two solutions of mean-De turck flow constructed above.
\end{prop}
\begin{pf} Let
$\psi(\bar{z_1},\bar{z_2})=d^{2}_{\bar{M}}(\bar{z}_1,\bar{z}_2)$
be the square of the distance function on $\bar{M}$ which is
viewed as a function of $(\bar{z}_1,\bar{z}_2)\in \bar{M}\times
\bar{M}.$  Set
$u(y,t)=d^{2}_{\bar{M}}(\bar{X}_1(y,t),\bar{X}_2(y,t)).$  Let
$\Delta_{k}=g_{k}^{\alpha\beta} \hat{{\nabla}}_{\alpha}
\hat{\nabla}_{\beta}$ for $k=1,2$. By direct computation, we have
  \begin{eqnarray*}\frac{\partial}{\partial t}u(y,t)
  &= &2d_{\bar{M}}(\bar{X}_1,\bar{X}_2)\frac{\partial d}{\partial
  \bar{z_{1}}^{\bar{\xi}}} \triangle_1\bar{X}_{1}^{\bar{\xi}}
  +2d_{\bar{M}}(\bar{X}_1,\bar{X}_2)\frac{\partial d}{\partial
  \bar{z_{2}}^{\bar{\zeta}}} \triangle_2\bar{X}_{2}^{\bar{\zeta}},\\
  g_{1}^{\alpha\beta} \hat{{\nabla}}_{\alpha} \hat{\nabla}_{\beta}u(y,t)
  &=& 2 d_{\bar{M}}(\bar{X}_1,\bar{X}_2) [\frac{\partial d}{\partial
  \bar{z_{1}}^{\bar{\xi}}} \triangle_1\bar{X}_{1}^{\bar{\xi}}+\frac{\partial d}{\partial
  \bar{z_{2}}^{\bar{\zeta}}}
  \triangle_1\bar{X}_{2}^{\bar{\zeta}}]+Hess
  (\psi)(Z_{\alpha},Z_{\beta})g_{1}^{\alpha\beta},
  \end{eqnarray*}
  where $Z_{\alpha}=\frac{\partial \bar{X}_{1}^{\bar{\xi}}}
  {\partial y^{\alpha}}\frac{\partial}{\partial \bar{z}_{1}^{\bar{\xi}}}
  +\frac{\partial \bar{X}_{2}^{\bar{\zeta}}}
  {\partial y^{\alpha}}\frac{\partial}{\partial
  \bar{z}_{2}^{\bar{\zeta}}} \in T_{(\bar{X}_{1},\bar{X}_{2})}\bar{M}\times
  \bar{M},$ $\alpha=1,2\cdots, n$ are vector fields on $\bar{M}\times \bar{M}.$
Combining these two formulas, we have $$[\frac{\partial}{\partial
t}-g_{1}^{\alpha\beta} \hat{{\nabla}}_{\alpha}
\hat{\nabla}_{\beta}]u(y,t)
  = -2d_{\bar{M}}(\bar{X}_1,\bar{X}_2)\frac{\partial d}{\partial
  \bar{z}_{2}^{\bar{\zeta}}} ((\triangle_1-\triangle_2)\bar{X}_{2})^{\bar{\zeta}}-Hess
  (\psi)(Z_{\alpha},Z_{\beta})g_{1}^{\alpha\beta}. \eqno{(5.8)}
  $$

  Note that
    \begin{equation*}\tag{5.9}
    \begin{split}
(\triangle_1-\triangle_2)\bar{X}_{2}&=g_{1}^{\alpha\beta}
\hat{{\nabla}}_{\alpha}
\hat{\nabla}_{\beta}\bar{X}_2-g_{2}^{\alpha\beta}
\hat{{\nabla}}_{\alpha}\hat{\nabla}_{\beta}\bar{X}_2\\[2mm]
&=g_{1}^{\alpha\gamma}g_{2}^{\beta\delta}({g_{2}}_{\delta\gamma}-{g_{1}}_{\delta\gamma})
\hat{{\nabla}}_{\alpha}
\hat{\nabla}_{\beta}\bar{X}_2,\\[2mm]
\hat{{\nabla}}_{\alpha} \hat{\nabla}_{\beta}\bar{X}_2&=
{\nabla_{2}}_{\alpha}{\nabla_{2}}_{\beta}\bar{X}_{2}+
(\hat{\Gamma}-\Gamma_{2})\ast \nabla\bar{X}_{2},
\end{split}
\end{equation*}
where $\Gamma_2$ and $\nabla_2$ are the christoffel symbol and the
covariant derivative of the metric ${g_{2}}_{\alpha\beta}(y,t).$

For each $y\in M^{n}$ and $t\in[0,T]$, if
$\bar{X}_{1}(y,t)\neq\bar{X}_{2}(y,t)$, denote the minimal
geodesic on $\bar{M}$ from $\bar{X}_{1}(y,t)$ to
$\bar{X}_{2}(y,t)$ by $\sigma$, and denote the parallel
translation of $\bar{M}$ along $\sigma$ by $P_{\sigma}$, then we
have
\begin{equation*}\tag{5.10}
\begin{split}
{g_{1}}_{\delta\gamma}(y,t)-{g_{2}}_{\delta\gamma}(y,t)&=
\langle\bar{X_{1}}_{\ast}(\frac{\partial}{\partial
y^{\delta}}),\bar{X_{1}}_{\ast}(\frac{\partial}{\partial
y^{\gamma}})\rangle_{\bar{g}}-\langle\bar{X_{2}}_{\ast}(\frac{\partial}{\partial
y^{\delta}}),\bar{X_{2}}_{\ast}(\frac{\partial}{\partial
y^{\gamma}})\rangle_{\bar{g}}\\
&=\langle\bar{X_{1}}_{\ast}(\frac{\partial}{\partial
y^{\delta}}),\bar{X_{1}}_{\ast}(\frac{\partial}{\partial
y^{\gamma}})\rangle_{\bar{g}}-\langle
P_{\sigma}^{-1}(\bar{X_{2}}_{\ast}(\frac{\partial}{\partial
y^{\delta}})),P_{\sigma}^{-1}(\bar{X_{2}}_{\ast}(\frac{\partial}{\partial
y^{\gamma}}))\rangle_{\bar{g}}\\
&=\langle\bar{X_{1}}_{\ast}(\frac{\partial}{\partial
y^{\delta}})-P_{\sigma}^{-1}(\bar{X_{2}}_{\ast}(\frac{\partial}{\partial
y^{\delta}})),\bar{X_{1}}_{\ast}(\frac{\partial}{\partial
y^{\gamma}})\rangle_{\bar{g}}\\&\ \ +\langle
P_{\sigma}^{-1}(\bar{X_{2}}_{\ast}(\frac{\partial}{\partial
y^{\delta}})),\bar{X_{1}}_{\ast}(\frac{\partial}{\partial
y^{\gamma}})-P_{\sigma}^{-1}(\bar{X_{2}}_{\ast}(\frac{\partial}{\partial
y^{\gamma}}))\rangle_{\bar{g}}.
\end{split}
\end{equation*}
If $\bar{X}_{1}(y,t)=\bar{X}_{2}(y,t),$ $P_{\sigma}=Identity$, the
above formula still holds.

In the following argument, we compute norms by using the metrics
$g_{1}$ and $\bar{g}.$  For example
$$ |\hat{\Gamma}-\Gamma_{2}|^{2}=
(\hat{\Gamma}-\Gamma_{2})^{\gamma}_{\alpha\beta}(\hat{\Gamma}-\Gamma_{2})^{\gamma'}_{\alpha'\beta'}{g_1}_{\gamma\gamma'}
{g_1}^{\alpha\alpha'}{g_1}^{\beta\beta'}$$ and
$$
|\nabla^{2}_{2}\bar{X}_{2}|^{2}=
\bar{g}_{\bar{\xi}\bar{\zeta}}g_{1}^{\alpha\alpha'}g_{1}^{\beta\beta'}
{\nabla_2}_{\alpha}{\nabla_2}_{\beta}\bar{X}_{2}^{\bar{\xi}}
{\nabla_2}_{\alpha'}{\nabla_2}_{\beta'}\bar{X}_{2}^{\bar{\zeta}}.
$$
We denote by $C$ various constants depending only on the constants
$\bar{C}$, $T$, $\bar{n}$ and $\bar{\delta}$ in the main theorem
1.1. Then by (5.3), we have
\begin{equation*}\tag{5.11}
\begin{split}
&|\hat{\Gamma}-\Gamma_{2}|\leq C,\\
&|\hat{\nabla}^{2}\bar{X}_{2}|\leq
C|\hat{\Gamma}-\Gamma_{2}|+C|{\nabla}^{2}_{2}\bar{X}_{2}|\leq C,\\
&|g_{2}|+|g_{2}^{-1}|\leq C,
\end{split}
\end{equation*}
where $|{\nabla}^{2}_{2}\bar{X}_{2}|$ is just the norm of the
second fundamental form of $X_2:M^{n}\rightarrow
\bar{M}^{\bar{n}}$ which is bounded by $\bar{C}$. Combining (5.9)
(5.10) and (5.11), we have
\begin{equation*}\tag{5.12}
\begin{split}
|(\triangle_1-\triangle_2)\bar{X}_{2}|^{2}&\leq
Cg_{1}^{\delta\gamma}\langle\bar{X_{1}}_{\ast}(\frac{\partial}{\partial
y^{\delta}})-P_{\sigma}^{-1}(\bar{X_{2}}_{\ast}(\frac{\partial}{\partial
y^{\delta}})),\bar{X_{1}}_{\ast}(\frac{\partial}{\partial
y^{\gamma}})-P_{\sigma}^{-1}(\bar{X_{2}}_{\ast}(\frac{\partial}{\partial
y^{\gamma}}))\rangle_{\bar{g}}.
\end{split}
\end{equation*}

  By choosing an orthonormal frame at $y$  so
 that ${g_{1}}_{\alpha\beta}=\delta_{\alpha\beta},$ then we have
$$Hess
  (\psi)(Z_{\alpha},Z_{\beta})g_{1}^{\alpha\beta}=\sum_{\alpha=1}^{n} Hess(\psi)(Z_{\alpha},Z_{\alpha}).$$
Note that
$$
Z_{\alpha}={Z_\alpha}_{1}+{Z_\alpha}_{2}, \ \ \ \ \ \ \mbox{for}\
\ \ \alpha=1,2,\cdots, n,
$$
where ${Z_\alpha}_{1}=\frac{\partial \bar{X}_{1}^{\bar{\xi}}}
  {\partial y^{\alpha}}\frac{\partial}{\partial
  \bar{z}_{1}^{\bar{\xi}}}= \bar{X_1}_{\ast}(\frac{\partial}{\partial
  y^{\alpha}})$ and ${Z_\alpha}_{2}=\frac{\partial \bar{X}_{2}^{\bar{\zeta}}}
  {\partial y^{\alpha}}\frac{\partial}{\partial
  \bar{z}_{2}^{\bar{\zeta}}}= \bar{X_2}_{\ast}(\frac{\partial}{\partial
  y^{\alpha}})$.

Recall that by Theorem 2.2 (ii), there is a constant $C$ such that
if $d_{\bar{M}}(\bar{z}_1,\bar{z}_2)\leq
  \min\{\frac{1}{4\sqrt{\bar{C}}},\frac{\bar{\delta}}{2}\}$, we have
  $$ (\nabla^{2}  d^{2})(Z,Z)\geq 2|Z_1-P_{\sigma}^{-1}Z_2|^{2}-C|Z|^{2}d^{2}\ \ \ \mbox{for all }
 Z\in T_{(\bar{z}_1,\bar{z}_2)}\bar{M}^{\bar{n}}\times \bar{M}^{\bar{n}},$$
 where $Z=Z_1+Z_2,$ $Z_1\in T_{\bar{z}_1}\bar{M}^{\bar{n}},$ $Z_2\in
 T_{\bar{z}_2}\bar{M}^{\bar{n}}.$ Hence if $d_{\bar{M}}(\bar{X}_1,\bar{X}_2)\leq
  \min\{\frac{1}{4\sqrt{\bar{C}}},\frac{\bar{\delta}}{2}\}$, then
  $$
  \sum_{\alpha=1}^{n}Hess(\psi)(Z_{\alpha},Z_{\alpha})\geq  \sum_{\alpha=1}^{n}2|\bar{X_1}_{\ast}(\frac{\partial}{\partial y^{\alpha}})-
  P_{\sigma}^{-1}\bar{X_2}_{\ast}(\frac{\partial}{\partial
  y^{\alpha}})|^{2}-Cd_{\bar{M}}(\bar{X}_{1},\bar{X}_{2})^{2}
  \eqno{(5.13)}
  $$
  since $|Z_{\alpha}|\leq C.$

  Combining (5.8), (5.12) and (5.13),  if  $u^{\frac{1}{2}}\leq
  \min\{\frac{1}{4\sqrt{\bar{C}}},\frac{\bar{\delta}}{2}\},$ then we have
  \begin{equation*} \tag{5.14}
  \begin{split}
(\frac{\partial}{\partial t}-g_{1}^{\alpha\beta}
\hat{{\nabla}_{\alpha}} \hat{\nabla}_{\beta})u(y,t)
  &\leq Cd_{\bar{M}}(\bar{X}_{1},\bar{X}_{2})\sum_{\alpha=1}^{n}2|\bar{X_1}_{\ast}(\frac{\partial}{\partial y^{\alpha}})-
  P_{\sigma}^{-1}\bar{X_2}_{\ast}(\frac{\partial}{\partial
  y^{\alpha}})|\\&\ \ \ \ -2\sum_{\alpha=1}^{n}|\bar{X_1}_{\ast}(\frac{\partial}{\partial y^{\alpha}})-
  P_{\sigma}^{-1}\bar{X_2}_{\ast}(\frac{\partial}{\partial
  y^{\alpha}})|^{2}+Cd_{\bar{M}}^{2}(\bar{X}_{1},\bar{X}_{2})\\
  &\leq Cu.\end{split}
  \end{equation*}

Now we show that $u^{\frac{1}{2}}\leq
  \min\{\frac{1}{4\sqrt{\bar{C}}},\frac{\bar{\delta}}{2}\}$ on
  some time interval $[0,T_5]$.

  For any $(y,t)\in \hat{M}\times[0,T_4]$, we have
\begin{equation*} \tag{5.15}
  \begin{split}
u^{\frac{1}{2}}(y,t)
  \leq & d_{\bar{M}}({X}_{1}\circ F_1^{-1}(y,t),{X}_{1}\circ
  F_1^{-1}(y,0))+d_{\bar{M}}({X}_{1}\circ F_1^{-1}(y,0),{X}_{2}\circ
  F_2^{-1}(y,0))\\[2mm]
  & +d_{\bar{M}}({X}_{2}\circ F_2^{-1}(y,t),{X}_{2}\circ
  F_2^{-1}(y,0))\\[2mm]
  \triangleq & I_1+I_2+I_3.\end{split}
  \end{equation*}
By the mean curvature flow equation (1.1), we know
\begin{equation*}
  \begin{split}
I_2  \leq
d_{\bar{M}}({X}_{1}(y,t),{X}_{1}(y,0))+d_{\bar{M}}({X}_{2}(y,t),{X}_{2}(y,0))\leq
 2\sqrt{n}\bar{C}t.\end{split}
  \end{equation*}
  By (4.4), (4.23),  for any $x_1, x_2\in M^{n}$, we get
  $$
  \frac{\partial}{\partial t}\hat{d}(F_1(x_1,t),F_1(x_2,t))\geq -
  C,
  $$
this implies
$$
\hat{d}(x_1,x_2)\leq \hat{d}(F_1(x_1,t),F_1(x_2,t))+Ct.
\eqno{(5.16)}
$$
By (5.16) and Lemma 4.2, it follows
\begin{equation*}
  \begin{split}
I_1  = &d_{\bar{M}}({X}_{1}\circ F_1^{-1}(y,t),{X}_{1}\circ
F_1^{-1}(y,0))\\[2mm]
\leq & d_{(M,g^{1}(\cdot,t))}( F_1^{-1}(y,t), y)\\[2mm]
\leq & C \hat{d}( F_1^{-1}(y,t), y) \\[2mm] \leq & Ct+
C\hat{d}(y,F_1(y,t))\\[2mm]
 \leq & C\sqrt{t}.\end{split}
  \end{equation*}
  The estimate of $I_3$ is similar.
Therefore, we have
$$
u^{\frac{1}{2}}(y,t)\leq C\sqrt{t} \eqno{(5.17)}
$$
for some constant $C$ depending only on $\bar{C}$, $\bar{\delta}$,
$T$ and $\bar{n}$.

Although $g_{1}^{\alpha\beta} \hat{{\nabla}_{\alpha}}
\hat{\nabla}_{\beta}$ is not the standard Laplacian, the maximum
principle is still applicable.  For completeness, we include the
proof in the following.

 Since the curvature of $(M,\hat{g})$ is bounded,
it is well-known that there is a function $\varphi$ such that
\begin{equation*}
\begin{split}
\frac{1}{C}(1+d_{\hat{g}}(y_0,y))&\leq \varphi(y)\leq C
(1+d_{\hat{g}}(y_0,y))\\
|\hat{\nabla}\varphi|+|\hat{\nabla}^{2}\varphi|&\leq C.
\end{split}
\end{equation*}
Note $g_{1}$ is equivalent to $\hat{g}$. For any small
$\varepsilon>0$ and big $A>0$,  we have
 \begin{equation*}
  \begin{split}
(\frac{\partial}{\partial t}-g_{1}^{\alpha\beta}
\hat{{\nabla}_{\alpha}}
\hat{\nabla}_{\beta})(e^{-Ct}u(y,t)-\varepsilon e^{At}\varphi)\leq-
\frac{\varepsilon A}{2} e^{At}\varphi<0.
\end{split}
  \end{equation*}
Then the classical maximum principle implies that for any fixed
$t_0$ the maximal value of $(e^{-Ct}u(y,t)-\varepsilon
e^{At}\varphi)$ on $M\times [0,t_0]$ can not be achieved for any
point $(y,t)$ with $0<t\leq t_0.$ Hence $e^{-Ct}u(y,t)-\varepsilon
e^{At}\varphi\leq 0$ for any $t\in [0,T_5]$ for some $T_5>0.$  Let
$\varepsilon\rightarrow 0,$ we conclude that
 $ u\equiv 0
 $
 on $[0,T_5].$
This implies $\bar{X}_{1}=\bar{X}_2,$ on $ M\times[0,T_5].$ We
complete the proof of Proposition 5.1. \end{pf}

%%%%%%%%%%%%%%%%%%%%%%%%%%%%%%%%%%%%%%%%%%%%%%%%%%%%%%%%%%%%%%%%%%%%%%%%%%%%%%%%%%%%%%%
\section{Proof of the uniqueness theorem 1.1}
\qquad  Now we are ready to prove Theorem 1.1. Let $X_1(x,t)$ and
$X_2(x,t)$ be two solutions of MCF with bounded second fundamental
form and with the same initial data. We solve the corresponding
harmonic map flow (5.4) (5.5)(with the same target $(M,\hat{g})$,
$\hat{g}=g_1(T)$) respectively to obtained two solutions
$F_1(x,t)$ and $F_2(x,t)$ on some common time interval. Then
$\bar{X}_1=X_1\circ F_1^{-1}$ and $\bar{X}_2=X_2\circ F_2^{-1}$
are two solutions to the mean-De Turck flow with the same initial
value. By Proposition 5.1 we know $\bar{X}_1\equiv\bar{X}_2$ on
$[0,T_5]$. So in order to prove $X_1(x,t)\equiv X_2(x,t)$, we only
need to show $F_1\equiv F_2$.

We know
\begin{eqnarray*}
 \Delta_1 F_1^{\alpha}&=&g_1^{\beta\gamma}(\hat{\Gamma}^{\alpha}_{\beta\gamma}-{\Gamma}^{\alpha}_{1\beta\gamma})\circ F_1,\\
 \Delta_2 F_2^{\alpha}&=&g_2^{\beta\gamma}(\hat{\Gamma}^{\alpha}_{\beta\gamma}-{\Gamma}^{\alpha}_{2\beta\gamma})\circ F_2.
\end{eqnarray*}
Since $\bar{X}_1\equiv\bar{X}_2,$ we know
$g_{1\alpha\beta}(y,t)=g_{2\alpha\beta}(y,t)$ on $[0,T_5]$, and
the
 vector fields $V_1\equiv V_2$ on the target, where
\begin{eqnarray*}
 V_1^{\alpha}&=&g_1^{\beta\gamma}(\hat{\Gamma}^{\alpha}_{\beta\gamma}-{\Gamma}^{\alpha}_{1\beta\gamma}),\\
 V_2^{\alpha}&=&g_2^{\beta\gamma}(\hat{\Gamma}^{\alpha}_{\beta\gamma}-{\Gamma}^{\alpha}_{2\beta\gamma}).
\end{eqnarray*}
Therefore, the two families of maps $F_1$ and $F_2$ satisfy
 the same ODE with the same initial value:
  \begin{equation*}
  \left\{
  \begin{split}
 \frac{\partial}{\partial t}F_1&=V\circ F_1\\
 F_1(\cdot,0)&=Identity,
 \end{split}
 \right.
 \end{equation*} and
\begin{equation*}
  \left\{
  \begin{split}
 \frac{\partial}{\partial t}F_2&=V\circ F_2\\
 F_2(\cdot,0)&=Identity.
\end{split}
 \right.
 \end{equation*}
 So  for any $x\in M^n$, letting $\gamma$ be a shortest geodesic( parametrized by arc length)
 on the target with $\gamma(0)=F_1(x,t)$ and $\gamma(l)=F_2(x,t)$, we have
 \begin{eqnarray*}
 \frac{\partial}{\partial t}\hat{d}(F_1(x,t),F_2(x,t))&=&\langle
 V,{\gamma}^{\prime}(l)\rangle-\langle
 V,{\gamma}^{\prime}(0)\rangle\\
 &=&\langle
 P_{\gamma}^{-1}V-V,{\gamma}^{\prime}(0)\rangle\\
 &\leq& \sup_{y\in\gamma}|\hat{\nabla}V|(y,t)
 \hat{d}(F_1(x,t),F_2(x,t)),
 \end{eqnarray*}
 where $P_{\gamma}^{-1}V$ is the parallel transport of
 $V(F_2(x,t),t)$ along the geodesic $\gamma$ back to the
 tangent space of the point $F_1(x,t)$. We have seen in the proof of Corollary 4.6 that $ \sup\limits_{y\in\gamma}|
  \hat{\nabla}V|(y,t)$ $\leq \frac{C}{\sqrt{t}}$ for some $C$ depending on $x$ but independent of $t$. Since $\hat{d}(F_1(x,0),F_2(x,0))\equiv 0$,
 we conclude that
 $$
 F_1(x,t)\equiv F_2(x,t).
 $$
 So we have proved $X_1(x,t)=X_2(x,t)$, for all $x\in M$ and
 $t\in [0,T_{5}]$. Clearly, we can extend the
 interval $[0,T_{5}]$ to the whole $[0,T]$ by applying the same argument on $[T_{5},T]$.

 The proof of Theorem 1.1 is completed. $\hfill\Box$

Corollary 1.2 is a direct consequence of Theorem 1.1. Indeed, let
$\bar{\sigma}$ and  $\sigma$ be two isometries of
$(\bar{M}^{\bar{n}},\bar{g})$ and $({M}^{{n}},{g})$ respectively
such that $(\bar{\sigma}\circ X_0)(x)=(X_0\circ \sigma)(x) $ for
any $x\in M^{n}.$  Since $\bar{\sigma}\circ X_t$ and $ X_t\circ
\sigma$ are two solutions to the MCF (1.1) with bounded second
fundamental form on $M^{n}\times [0,T]$ and with the same initial
value, then by Theorem 1.1, we have
$$(\bar{\sigma}\circ X_t)(x)=(X_t\circ \sigma)(x)
$$
for any $x\in M^{n}$ and $t\in[0,T].$ The proof of the Corollary
1.2 is completed.  $\hfill\Box$

\section{Pseudolocality Theorem}

\qquad We begin with a few terminologies for the sake of
convenience. An $n$-dimensional submanifold $M\subset\bar{M}$ is
said to be a local $\delta$- Lipschitz graph of radius $r_0$ at $P
\in {M}$, if there is a normal coordinate system $(y^{1}\cdots
y^{\bar{n}})$ of $\bar{M}$ around $P$ with
$T_{P}M=$span$\{\frac{\partial}{\partial
y^{1}},\cdots,\frac{\partial}{\partial y^{n}}\}$, a vector valued
function $F:
\{y'=(y^1,\cdots,y^n)\mid(y^1)^{2}+\cdots+(y^n)^{2}<r_0^{2}\}\rightarrow
\mathbb{R}^{\bar{n}-n} $ with $F(0)=0$, $|DF|(0)=0$ such that
 $M\cap\{|y'|<r_0\}=\{(y',F(y'))\mid |y'|<r_0\}$ and
$|DF|^{2}(y')=\sum_{i,\beta}\frac{\partial F^{\beta}}{\partial
y^{i}}\frac{\partial F^{\beta}}{\partial y^{i}}< \delta^{2}$. The
submanifold $M_0$ is said to be graphic in the ball
$B_{\bar{M}}(x_0,r_0)$, if the above holds for $\delta=\infty.$

We say a submanifold $M\subset \bar{M}$ is properly embedded in a
ball $B_{\bar{M}}(x_0, r_0)$ if either $M$ is closed or $\partial
M$ has distance $\geq r_0$ from $x_0.$ We say a submanifold
$M\subset \bar{M}$ is properly embedded in $\bar{M}$ if either $M$
is closed or there is an $x_0\in \bar{M}$ such that $M$ is
properly embedded in $B_{\bar{M}}(x_0, r_0)$ for any $r_0>0.$ It
is clear that if $\bar{M}$ is complete and $M$ is properly
embedded in  $\bar{M},$ then $M$ is complete. A properly embedded
submanifold $M$ is said to be uniform graphic with radius $r_0$ if
for any $x_0\in M$ it is graphic in the ball
$B_{\bar{M}}(x_0,r_0).$

The following lemma says that if the second fundamental form is
controlled, then  (a piece of) the sub-manifold is a local
$\delta$-Lipschitz graph of suitable radius.
\begin{lem}
 Let $\bar{M}$ be an $\bar{n}-$dimensional complete  Riemannian manifold
  satisfying  $$|\bar{R}m|+|\bar{\nabla}\bar{R}m|(x)\leq
 \bar{C}, \ \ \ \ \ \  inj(\bar{M})\geq i_{0}>0. \ \ \ \
 $$
 There exists a constant $C_1>0$ with the following property. Let $\{x^{1},\cdots,x^{\bar{n}}\}$ be normal
   coordinates of $\bar{M}$ of radius $r_0$ around $x_0$ with $T_{x_0}M=
   span\{\frac{\partial}{\partial x^{1}},\cdots,\frac{\partial}{\partial
   x^{n}}\},$ where $M$ is an $n$-dimensional submanifold properly embedded in $B_{\bar{M}}(x_0,r_0),$
    $x_0\in M,$ $r_0\leq\frac{1}{C_1},$ and the second fundamental form $|A|\leq
   \frac{1}{r_0}.$ Then there exists a map
   $F:\{(x^{1},\cdots,x^{{n}})\mid{ (x^{1}}^{2}+\cdots+{x^{n}}^{2})^{\frac{1}{2}}<\frac{r_0}{96}\}\rightarrow \mathbb{R}^{\bar{n}-n}$
   with  $F(0)=0,$ $|DF|(0)=0$ such that the connected component
   containing $x_0$
   of
   $M\cap \{(x^{1},\cdots,x^{\bar{n}})\mid
   ({x^{1}}^{2}+\cdots+{x^{n}}^{2})^{\frac{1}{2}}<\frac{r_0}{96}\}$ can be written as a graph $\{(x^{'},F(x'))\mid
   {|x'|}=({x^{1}}^{2}+\cdots+{x^{n}}^{2})^{\frac{1}{2}}<\frac{r_0}{96}\}$
   and
   $$|DF|(x')\leq \frac{36}{r_0}|x'|,\eqno{(7.1)}$$ $x'=(x^{1},\cdots,x^{n})\in B_{\mathbb{R}^{n}}(0,\frac{r_0}{96}),$
   where $|DF|(x')^{2}=
   \sum_{i=1}^{n}\sum_{\alpha=n+1}^{\bar{n}}\frac{\partial
   F^{\alpha}}{\partial x^{i}}\frac{\partial F^{\alpha}}{\partial
   x^{i}}(x').$
   \end{lem}
\begin{pf} Let
$X=(X^1,\cdots,X^{\bar{n}})=(x',F(x')),x'=(x^1,\cdots,x^{n}),$ be
a graph representation of  the local isometric embedding of the
connected component containing $x_0$ of $M\cap
\{(x^{1},\cdots,x^{\bar{n}})\mid
   ({x^{1}}^{2}+\cdots+{x^{n}}^{2})^{\frac{1}{2}}<r_1\}$(for some $r_1\leq\frac{r_0}{96}$)
    into $\bar{M}$ under the exponential map.

Define $$ |\nabla
F|^{2}=\sum_{i,j=1}^{n}\sum_{\alpha=n+1}^{\bar{n}}\frac{\partial
F^{\alpha}}{\partial x^{i}}\frac{\partial F^{\alpha}}{\partial
x^{j}}g^{ij},|D
F|^{2}=\sum_{i=1}^{n}\sum_{\alpha=n+1}^{\bar{n}}\frac{\partial
F^{\alpha}}{\partial x^{i}}\frac{\partial F^{\alpha}}{\partial
x^{i}}.
$$
By choosing $C_1$ large, we have
$$
\frac{1}{2}\delta_{\alpha\beta}\leq \bar{g}_{\alpha\beta}\leq 2
\delta_{\alpha\beta}, \ \
|\bar{\Gamma}^{\gamma}_{\alpha\beta}|\leq 1,\ \
\frac{1}{2}\delta_{ij}\leq g_{ij}\leq 2(1+|DF|^{2}) \delta_{ij}.
$$
For $\alpha\geq n+1,$ $i,j\leq n,$ recall the coefficients of the
second fundamental form is given by
$$
A^{\alpha}_{ij}=\frac{\partial X^{\alpha}}{\partial x^{i}\partial
x^{j}}-\Gamma^{k}_{ij}\frac{\partial X^{\alpha}}{\partial
x^{k}}+\bar{\Gamma}^{\alpha}_{\beta\gamma}\frac{\partial
X^{\beta}}{\partial x^{i}}\frac{\partial X^{\gamma}}{\partial
x^{j}}=\nabla^{2}_{ij}
F^{\alpha}+\bar{\Gamma}^{\alpha}_{\beta\gamma}\frac{\partial
X^{\beta}}{\partial x^{i}}\frac{\partial X^{\gamma}}{\partial
x^{j}}.
$$
Note that $$|\bar{\Gamma}^{\alpha}_{\beta\gamma}\frac{\partial
X^{\beta}}{\partial x^{i}}\frac{\partial X^{\gamma}}{\partial
x^{j}}|^{2}=\bar{\Gamma}^{\alpha'}_{\beta'\gamma'}\frac{\partial
X^{\beta'}}{\partial x^{i}}\frac{\partial X^{\gamma'}}{\partial
x^{j}}\bar{\Gamma}^{\alpha}_{\beta\gamma}\frac{\partial
X^{\beta}}{\partial x^{k}}\frac{\partial X^{\gamma}}{\partial
x^{l}}g^{ik}g^{jl}\bar{g}_{\alpha\alpha'}\leq C(\bar{n}),$$

\begin{equation*}
\begin{split}
|\nabla ^{2}F|^{2}&=\sum_{\alpha,\beta\geq n+1;i,j,k,l\leq
n}\nabla^{2}_{ij}F^{\alpha}\nabla^{2}_{kl}F^{\beta}\delta_{\alpha\beta}g^{ik}g^{jl}\\
&\leq 4(|A|^{2}+C(\bar{n}))\\
&\leq 4r_{0}^{-2}+C(\bar{n}),
\end{split}
\end{equation*}
and $$|\nabla|\nabla F||\leq |\nabla ^{2}F|.
$$
This implies
$$|\nabla F|(\cdot)\leq 3r_{0}^{-1}d_{M}(x_0,\cdot). \eqno{(7.2)}
$$
 Since $g_{ij}\leq
2(\delta_{ij}+\frac{\partial F^{\alpha}}{\partial
x^{i}}\frac{\partial F^{\alpha}}{\partial x^{j}})\leq
2(1+|DF|^{2})\delta_{ij},$ it follows that
$$
|\nabla F|^{2}\geq \frac{1}{4}\frac{|D F|^{2}}{1+|D F|^{2}}
$$
and $$|D F|^{2}\leq \frac{4|\nabla F|^{2}}{1-4|\nabla
F|^{2}}.\eqno{(7.3)}$$ Combining (7.2)and (7.3), it follows that
$$|D F|(\cdot)\leq 9 r_0^{-1}d_{M}(x_0,\cdot) \ \ \ \ \ \text{on } \ \ B_{M}(x_0,\frac{r_0}{24}).$$

Since  $d_{M}(x_0,\cdot) \leq 2 d_{\bar{M}}(x_0,\cdot)$ by (2.5),
we have
\begin{equation*}
\begin{split}
|D F|(\cdot)\leq 18
r_{0}^{-1}\sup_{B_{M}(0,\frac{r_0}{24})}(1+|DF|)|x'|\leq 36
r_{0}^{-1}|x'|,
\end{split}
\end{equation*}
and we conclude that
$$
|D F|(x')\leq 36 r_{0}^{-1}|x'|, \ \ \ \ \ \text{whenever} \ \
|x'|\leq \frac{r_0}{96}.
$$
 The above argument shows that there is
$C_1>0$ such that under the exponential map,  once the connected
component of $M$ can be expressed as a graph $(x',F(x'))$ on
$B_{\mathbb{R}^{n}}(0,r_1),$ for $r_1\leq \frac{r_0}{96},$ then
the estimate (7.1) holds. Hence the connected component of $M$ can
be expressed as a graph on the ball
$B_{\mathbb{R}^{n}}(0,\frac{r_0}{96}).$
\end{pf}

For future applications in pseudolocality theorem, we need a local
graph representation for mean curvature flow.

\begin{lem} Fix $k\geq 1.$ Let $\bar{M}$ be an $\bar{n}-$dimensional complete manifold
  satisfying  $$\sum\limits_{i=0}^{k+1}|\bar{\nabla}^{i}\bar{R}m|(x)\leq
 \bar{C}, \ \ \ \ \ \  inj(\bar{M})\geq i_{0}>0. \ \ \ \
 $$
 There exists a constant $C_1>0$ with the following property. Suppose
$M_s,$ $s\in [-r_0^{2},0]$ is a solution of MCF
  properly embedded in $B_{\bar{M}}(x_0,r_0),$
  $x_0\in M_0,$ $r_0\leq\frac{1}{C_1},$ with $\sum\limits_{i=0}^{k}|\nabla^{i} A|r_0^{i+1}\leq
1$ on   $B_{\bar{M}}(x_0,r_0).$ Denote by $x_0^{s}\in
   M_s$ the orbit of $x_0$.
   Let $\{x^{1},\cdots,x^{\bar{n}}\}$ be normal
  coordinates of $\bar{M}$ of radius $r_0$ around $x_0$ with $T_{x_0}M_0=
  span\{\frac{\partial}{\partial x^{1}},\cdots,\frac{\partial}{\partial
  x^{n}}\}$.
  Then there exist a family of smooth maps
  $F_s:\{(x^{1},\cdots,x^{{n}})\mid{ (x^{1}}^{2}+\cdots+{x^{n}}^{2})^{\frac{1}{2}}<\frac{r_0}{C_1}\}\rightarrow \mathbb{R}^{\bar{n}-n}$
  with  $F_0(0)=0,$ $|D_0F|(0)=0$, $\bar{exp}_{x_0}((0,F_s(0)))=x_0^{s}$
  such that
  the connected component
  of  $M_s\cap \{(x^{1},\cdots,x^{\bar{n}})\mid
  ({x^{1}}^{2}+\cdots+{x^{n}}^{2})^{\frac{1}{2}}<\frac{r_0}{C_1}\}$
   (under the exponential map $\bar{exp}_{x_0}$) containing $x_0^s$
  can be written as a graph $\{(x^{'},F_s(x'))\mid
  {|x'|}=({x^{1}}^{2}+\cdots+{x^{n}}^{2})^{\frac{1}{2}}<\frac{r_0}{C_1}\}$;
  moreover  we have $\sum\limits_{i=1}^{k+2}r_0^{i+1}|D^{i}F_s|\leq C_1.$
\end{lem}

Actually, by the MCF equation $\frac{\partial}{\partial
s}X=\triangle
   X$,   where $X=(x',F_s(x'))$ is the graph
   representations on $B(0,r_1)$ for some $r_1<\frac{r_0}{C_1},$
     we
   have  information on $|\frac{\partial}{\partial s}
   F_s|r_0+|\frac{\partial}{\partial s}DF_s|r_0^{2}\leq C_1$.
   This gives $|F_s(0)|\leq Csr_0^{-1}$ and $|DF_s|(0)\leq Csr_0^{-2}.$ Similarly,
   by integrating $|\nabla|\nabla F||\leq |\nabla ^{2}F|$,
   we know the graph representation holds in a ball of uniform
   radius $\frac{r_1}{C_1}.$ The higher derivative $D^{i}F$ can be estimated by $\sum_{j\leq
   i}|\nabla^{j}F|$ by definitions.
    $\hfill\Box$

Now we state the pseudolocality theorem for the MCF.
\begin{thm} Let $\bar{M}$ be an $\bar{n}$-dimensional complete manifold
satisfying $\sum\limits_{i=0}^{3}|\bar{\nabla}^{i}\bar{R}m|\leq
c_0^{2}$ and $inj(\bar{M})\geq i_0>0$. Then for every $\alpha>0$
there exist $\varepsilon>0$, $\delta>0$ with the following property.
Suppose we have a smooth solution to the mean curvature flow
$M_t\subset\bar{M}$ properly embedded in $B_{\bar{M}}(x_0, r_0)$ for
$t\in [0,T]$ with $0<T\leq \varepsilon^{2} r_0^{2}$, and assume that
at time zero, $M_0$ is a local $\delta$- Lipschitz graph of radius
$r_0$ at $x_0\in {M}_0$ with $r_0\leq\frac{i_0}{2}$. Then we have an
estimate of the second fundamental form
$$
|A|(x,t)^{2}\leq\frac{\alpha}{t}+(\varepsilon r_0)^{-2} \eqno{(7.4)}
$$
on $B_{\bar{M}}(x_0,\varepsilon r_0)\cap M_t$, for any $t\in
[0,T]$.
\end{thm}
\begin{pf} We argue by contradiction. By scaling we may assume $r_0=1$.
Suppose there exist fixed $c_0>0$, $i_0>0$, $\alpha>0$,
   and a sequence of $\varepsilon,\delta\rightarrow0$
   and smooth solutions to the mean curvature flow $M_t\subset\bar{M}$
   for $t\in [0,T]\subseteq[0,\varepsilon^{2}]$ such that at time zero, $M_0$ is a local $\delta$-
Lipschitz graph of radius $1$ at $x_0\in {M}$. But there is some
$(x_1,t_1)$ satisfying $0\leq t_1\leq T$ and $x_1\in
B_{\bar{M}}(x_0,\varepsilon)$ such that
$$
|A|(x_1,t_1)^{2}>\frac{\alpha}{t_1}+\varepsilon ^{-2}.
$$
Denote by $E_\alpha$ the set of points ($x,t$) satisfying
$|A|(x,t)^{2}\geq\frac{\alpha}{t}.$ Now we use the Perelman's
point-picking technique \cite{P1} to choose another point which
controls nearby points in its scale.
\begin{lem} For any $K>0$ with $K\varepsilon<\frac{1}{100n}$, let
$M_t$ be assumed as in the theorem, suppose $|A|(x_1,t_1)^{2}\geq
\frac{\alpha}{t_1}+\varepsilon^{-2}$ for some $(x_1,t_1)$
satisfying $0\leq t_1\leq T \leq \varepsilon^{2}$ and $x_1\in
B_{\bar{M}}(x_0,\varepsilon)$, then one can find
$(\bar{x},\bar{t})\in E_\alpha$ with $0< \bar{t}\leq T$,
$d_{\bar{M}}(x_0,\bar{x})\leq (2K+1)\varepsilon$ such that
$$
|A|(x,t)\leq 4 Q \eqno{(7.5)}
$$
whenever $\bar{t}-\frac{3}{4}\alpha Q^{-2}\leq t\leq \bar{t}$,
$d_{\bar{M}}(x,\bar{x})\leq KQ^{-1}$, where
$Q=|A|(\bar{x},\bar{t}).$
\end{lem}
 Firstly, we claim that there exists $(\bar{x},\bar{t})\in
E_\alpha$ with $0< \bar{t}\leq T$, $d_{\bar{M}}(x_0,\bar{x})\leq
(2K+1)\varepsilon$ such that$$ |A|(x,t)\leq 4 |A|(\bar{x},\bar{t})
$$
whenever $(x,t)\in E_{\alpha}$, $0\leq t\leq \bar{t}$,
$d_{\bar{M}}(x_0,{x})\leq
d_{\bar{M}}(x_0,\bar{x})+K|A|(\bar{x},\bar{t})^{-1}$.

The  argument is  by contradiction. If $(x_1,t_1)$ can not be
chosen for $(\bar{x},\bar{t})$, one can find $(x_2,t_2)\in
E_{\alpha}$ with $0\leq t_2\leq t_1$, $d_{\bar{M}}(x_0,{x_2})\leq
d_{\bar{M}}(x_0,{x_1})+K|A|({x_1},{t_1})^{-1}$,$|A|(x_2,t_2)>
4|A|(x_1,t_1)$. Inductively, we have a sequence of $(x_k,t_k)\in
E_{\alpha}$ with $0\leq t_k\leq t_{k-1}$,
$d_{\bar{M}}(x_0,{x_k})\leq
d_{\bar{M}}(x_0,{x_{k-1}})+K|A|({x_{k-1}}$, ${t_{k-1}})^{-1}$,
 $|A|(x_k,t_k)> 4|A|(x_{k-1},t_{k-1})$. Therefore we have
 $$
|A|(x_k,t_k)> 4^{k-1}|A|(x_{1},t_{1})\geq 4^{k-1}\varepsilon^{-1}
 $$
 and
$d_{\bar{M}}(x_0,{x_k})\leq
d_{\bar{M}}(x_0,{x_{1}})+K\sum_{i=1}^{\infty}(4^{i-1}|A|(x_{1},t_{1}))^{-1}\leq
(2K+1)\varepsilon <\frac{1}{2}$. Since the solution is smooth, we
get a contradiction as $k$ large enough.

For the chosen $(\bar{x},\bar{t})$,  if $(x,t)\notin E_{\alpha}$,
$\bar{t}-\frac{3}{4}\alpha Q^{-2}\leq t\leq \bar{t}$, then
$$
|A|^{2}(x,t)\leq \frac{\alpha}{t}\leq
\frac{\alpha}{\bar{t}-\frac{3}{4}\alpha Q^{-2}}\leq
 4Q^{2}.
$$
If $(x,t)\in E_{\alpha}$ and $d_{\bar{M}}(x,\bar{x})\leq
K|A|(\bar{x},\bar{t})^{-1}$, by above claim we still have the
estimate. The lemma is proved.

 Continuing the proof of  Theorem 7.3.

Choose $K=\frac{1}{\sqrt{\varepsilon}}$. Let $(\bar{x},\bar{t})$ be
the point obtained in Lemma 7.4. Consider the auxiliary functions
$$
\varphi(x,t)= (4\pi(\bar{t}-t))^{-\frac{n}{2}}
e^{-(1+\frac{1}{\varepsilon}(t-\bar{t}))\frac{d^{2}_{\bar{M}}(\bar{x},x)}{4(\bar{t}-t)}-\frac{n}{2\varepsilon}t},
\psi(x,t)=(1-\frac{d_{\bar{M}}(\bar{x},x)^{2}+3nt}{\rho^{2}})_{+}^{3}
$$
on $\bar{M}\times[0,\bar{t}]$, where
$\rho=\min\{\frac{1}{2},\frac{1}{c_0
\sqrt{e}},i_0,\sqrt{\varepsilon}\}$. They are also functions on $M$
by composing the inclusion maps.  We will compute their equations on
$M$. Since the sectional curvature of $\bar{M}$ satisfies
$-c_0^{2}\leq sec\leq c_0^{2}$, by comparison theorem and mean
curvature flow equation, we have
\begin{eqnarray*}
(\frac{\partial}{\partial
t}+\triangle)d_{\bar{M}}(\bar{x},\cdot)^{2}&=&4d_{\bar{M}}
\bar{\nabla} d_{\bar{M}}\cdot
{H}+tr(Hess(d^{2}_{\bar{M}}(\bar{x},\cdot))\mid_{TM})\\
&\geq &4d_{\bar{M}} \bar{\nabla} d_{\bar{M}}\cdot {H}+2n\frac{c_0 d_{\bar{M}}(\bar{x},\cdot)}{\tan c_0 d_{\bar{M}}(\bar{x},\cdot)}\\
&\geq &4d_{\bar{M}} \bar{\nabla} d_{\bar{M}}\cdot
{H}+2n(1-\frac{1}{2}c_0^{2}d^{2}_{\bar{M}}(\bar{x},\cdot)),\\
(\frac{\partial}{\partial
t}-\triangle)d_{\bar{M}}(\bar{x},\cdot)^{2}&=&-tr(Hess(d^{2}_{\bar{M}}(\bar{x},\cdot))\mid_{TM})\\
&\geq &-2n c_0 d_{\bar{M}}(\bar{x},\cdot){\coth (c_0
d_{\bar{M}}(\bar{x},\cdot))}\geq -3n
\end{eqnarray*}
whenever
$d_{\bar{M}}(\bar{x},\cdot)^{2}<\min\{\frac{1}{c_0^{2}e},i_0^{2}\}$,
$t\in[0,\bar{t}]$. Hence we have
$$(\frac{\partial}{\partial
t}-\triangle)\psi \leq 0\eqno{(7.6)}$$ and
\begin{equation*}\tag{7.7}
\begin{split} (\frac{\partial}{\partial t}+\triangle-|H|^{2})\varphi
&=\varphi[\frac{n}{2(\bar{t}-t)}-\frac{1+\frac{1}{\varepsilon}(t-\bar{t})}{4(\bar{t}-t)}(\frac{\partial}{\partial
t}+\triangle)d_{\bar{M}}(\bar{x},\cdot)^{2}
-\frac{(1+\frac{1}{\varepsilon}(t-\bar{t}))d_{\bar{M}}(\bar{x},\cdot)^{2}}{4(\bar{t}-t)^{2}}
\\[2mm]
& \ \ \ + \frac{(1+\frac{1}{\varepsilon}(t-\bar{t}))^{2}|\nabla
d_{\bar{M}}(\bar{x},\cdot)^{2}|^{2}}{16(\bar{t}-t)^{2}}-\frac{\frac{1}{\varepsilon}d_{\bar{M}}
(\bar{x},\cdot)^{2}}{4(\bar{t}-t)}-\frac{n}{2\varepsilon}-|H|^{2}]\\[2mm]
&\leq\varphi[-\frac{1+\frac{1}{\varepsilon}(t-\bar{t})}{(\bar{t}-t)}d_{\bar{M}}
\bar{\nabla} d_{\bar{M}}\cdot
{H}-\frac{(1+\frac{1}{\varepsilon}(t-\bar{t}))d_{\bar{M}}(\bar{x},\cdot)^{2}}{4(\bar{t}-t)^{2}}\\[2mm]
& \ \ \ + \frac{(1+\frac{1}{\varepsilon}(t-\bar{t}))^{2}|\nabla
d_{\bar{M}}(\bar{x},\cdot)^{2}|^{2}}{16(\bar{t}-t)^{2}}
-\frac{[\frac{1}{\varepsilon}-(1+\frac{1}{\varepsilon}(t-\bar{t}))n
c_0^{2}]d_{\bar{M}}(\bar{x},\cdot)^{2}}{4(\bar{t}-t)}-|H|^{2}]
\\[2mm]
&\leq
-|H+(1+\frac{1}{\varepsilon}(t-\bar{t}))\frac{d_{\bar{M}}(\bar{x},\cdot)\bar{\nabla}^{\perp}
d_{\bar{M}}(\bar{x},\cdot) }{2(\bar{t}-t)}|^{2}\varphi
\end{split}
\end{equation*}
whenever $d_{\bar{M}}(\bar{x},\cdot)<\rho$, $t\in[0,\bar{t}]$. We
used $0<1+\frac{1}{\varepsilon}(t-\bar{t})\leq 1.$ In the above and
following argument, we regard the mean curvature flow $M_t$ is a
smooth family of $F_t: M\rightarrow \bar{M},$  $(\varphi \psi) \circ
F_t $  is a $C^{2}$ function on  $M\times [0,\bar{t}]$ with compact
support in $M$. So $\int_{M_t}\varphi\psi=\int_{M}\varphi\psi dv_t$
is a $C^{2}$ function in $t$. Combining (7.6) and (7.7), we get the
monotonicity formula
$$
\frac{d}{dt}\int_{M_t}\varphi\psi\leq
-\int_{M_t}|H+(1+\frac{1}{\varepsilon}(t-\bar{t}))\frac{d_{\bar{M}}(\bar{x},\cdot)\bar{\nabla}^{\perp}
d_{\bar{M}}(\bar{x},\cdot) }{2(\bar{t}-t)}|^{2}\varphi\psi
\eqno{(7.8)}
$$
on $[0,\bar{t}]$. This implies
\begin{equation*}\tag{7.9}
\begin{split}
& \int_{\bar{t}-\frac{1}{2}\alpha Q^{-2}}^{\bar{t}}[\int_{M_t}
|H+(1+\frac{1}{\varepsilon}(t-\bar{t}))\frac{d_{\bar{M}}(\bar{x},\cdot)\bar{\nabla}^{\perp}
d_{\bar{M}}(\bar{x},\cdot) }{2(\bar{t}-t)}|^{2}\varphi\psi]dt\\
& \leq \int_{M_{\bar{t}-\frac{1}{2}\alpha
Q^{-2}}}\varphi\psi-\int_{M_{\bar{t}}}\varphi\psi.
\end{split}
\end{equation*}
Since the solution is smooth and properly embedded, $\psi$ is
compactly supported, we have
$\lim_{t\rightarrow\bar{t}_{-}}\int_{M_{{t}}}\varphi\psi=e^{-\frac{n}
{2\varepsilon}\bar{t}}(1-\frac{3n\bar{t}}{\rho^{2}})^{3}.$
 Now we
claim that there is $\beta>0$ such that as
$\varepsilon,\delta\rightarrow0$, we have
$$
\int_{M_{\bar{t}-\frac{1}{2}\alpha Q^{-2}}}\varphi\psi\geq
(1+\beta)e^{-\frac{n}{2\varepsilon}\bar{t}}(1-\frac{3n\bar{t}}{\rho^{2}})^{3}.\eqno{(7.10)}
$$
We still argue by contradiction. Suppose not, then there is a
subsequence of $\varepsilon,\delta\rightarrow0$ and $$
\int_{\bar{t}-\frac{1}{2}\alpha
Q^{-2}}^{\bar{t}}[\int_{M_t}|H+(1+\frac{1}{\varepsilon}(t-\bar{t}))\frac{d_{\bar{M}}(\bar{x},\cdot)\bar{\nabla}^{\perp}
d_{\bar{M}}(\bar{x},\cdot) }{2(\bar{t}-t)}|^{2}\varphi\psi dv]dt\leq
\beta\rightarrow 0. \eqno{(7.11)}
$$
Parabolic scaling the solution around $(\bar{x},\bar{t})$ with the
factor $Q$ and shifting the $\bar{t}$ to 0 and $\bar{x}$ to the
origin $O$, i.e. let
$(\tilde{M},\tilde{g})=(\bar{M},Q^{2}\bar{g})$ be the new target
manifold, $\tilde{M}_{s}=M_{\bar{t}+Q^{-2}s}$,
$-\frac{3}{4}\alpha\leq s\leq 0$ be the new family of
submanifolds, which is still solution of MCF. By (7.5), the
normalized second fundamental form satisfies $|\tilde{A}|\leq 4$
on $B_{\tilde{M}}(\bar{x},K)$, $-\frac{3}{4}\alpha \leq s\leq 0$.
By Theorem 3.2, we have
$|\tilde{\nabla}\tilde{A}|+|\tilde{\nabla}^{2}\tilde{A}|\leq
Const.$ on $B_{\tilde{M}}(\bar{x},\frac{K}{2})$,
$-\frac{5}{8}\alpha \leq s\leq 0$. Note that $K\rightarrow
\infty$.

  Now we are going to consider the convergence of the
MCF on changing target manifolds. We clarify the meaning of the
convergence in the following.

 Denote  the orbit of
$\bar{x}$ under MCF by $\bar{x}^{s}\in \tilde{M}_s$ such that
$\bar{x}^{0}=\bar{x}.$ Note the injectivity radius of the new
target $(\tilde{M},\tilde{g})$ tends to infinity as
$\varepsilon\rightarrow 0.$  Let $\{x^{1},\cdots,x^{\bar{n}}\}$ be
normal
   coordinates of $\tilde{M}$ of radius $\gg 1$ around $\bar{x}$ with $T_{\bar{x}}\tilde{M}_0=
   span\{\frac{\partial}{\partial x^{1}},\cdots,\frac{\partial}{\partial
   x^{n}}\},$  and $\tilde{g}_{\alpha\beta}$ be the metric coefficients of $\tilde{M}$ in this coordinates.
By \cite{Ha3}, we have
$|\tilde{g}_{\alpha\beta}-\delta_{\alpha\beta}|(x)\leq
CQ^{-2}|x|^{2}$ and  $
|\partial\tilde{g}_{\alpha\beta}|+|\partial^{2}\tilde{g}_{\alpha\beta}|\leq
C.$ By Arzela-Ascoli theorem, after taking a subsequence of
$\varepsilon\rightarrow0$, $\tilde{g}_{\alpha\beta}$ tends to $
\delta_{\alpha\beta}$ in $C^{2-\gamma}$ topology for any
$0<\gamma<1$.

    By  Lemma 7.2, there exist a family of maps
   $F_s:\{(x^{1},\cdots,x^{{n}})\mid{ (x^{1}}^{2}+\cdots+{x^{n}}^{2})^{\frac{1}{2}}<1\}\rightarrow \mathbb{R}^{\bar{n}-n}$
   with  $F_0(0)=0,$ $|DF_0|(0)=0$, such that the connected component
   containing $\bar{x}^s$
   of
   $\tilde{M}_s\cap \{(x^{1},\cdots,x^{\bar{n}})\mid
   ({x^{1}}^{2}+\cdots+{x^{n}}^{2})^{\frac{1}{2}}<1\}$ can be written as a graph $\{(x^{'},F_s(x'))\mid
   {|x'|}=({x^{1}}^{2}+\cdots+{x^{n}}^{2})^{\frac{1}{2}}<1\}.$ Moreover, we can show
   $$\sum\limits_{i=1}^{4}|D^{i}F|+\sum\limits_{i=1}^{2}(|\frac{\partial^{i}}{\partial s^{i}}F|
   +|D^{i}\frac{\partial F}{\partial
   s}|)\leq C,$$
where $D$ and the norm are  the natural differential and
   norm in Euclidean ordinates of $N\subset \mathbb{R}^{n}$ and the garget $\mathbb{R}^{\bar{n}}.$
   By Arzela-Ascoli theorem, $F(x',s)$ will converge to
   $F^{\infty}(x',s)$ in the topology of $\mathcal{C}^{\frac{3}{2}}(\overline{B(0,\frac{1}{2})}\times
   [-\frac{5\alpha}{8},0],\mathbb{R}^{\bar{n}})\cap \mathcal{C}^{3}(\overline{B(0,\frac{1}{2})},\mathbb{R}^{\bar{n}}).$

    If we set $X=(x',F(x'))$ being the
   map from $N:=B(0,1)$ to $\tilde{M},$ then the MCF equation
   can be written as
   $$
   \frac{\partial X}{\partial s}=\triangle X,
   $$
where $\triangle$ is the harmonic Laplacian defined by using the
induced metric $X^{*}\tilde{g}$ and the target metric $\tilde{g}.$
Since $X^{*}\tilde{g} $ is defined by $DF$ and $\tilde{g},$ after
taking  a subsequence of $\varepsilon\rightarrow 0,$ we know
$X^{*}\tilde{g}$ converges in
$\mathcal{C}^{1-\gamma}(\overline{B(0,\frac{1}{2})}\times
   [-\frac{5\alpha}{8},0])$ topology.

Denote by $\hat{M_s}=\tilde{M}_{s}\cap \exp_{\bar{x}}\{|x'|<1\},$
and $\hat{M}=\cup_{s\in [-\frac{\alpha}{2},0]}\hat{M_s}.$ By
summing up the above discussion, the piece  $\hat{M}$ of
$\tilde{M}$ containing $(\bar{x},0)$ will converge to a solution
of the MCF (in the classical sense) which is embedded on the
Euclidean space $\mathbb{R}^{\bar{n}}$ with
$|\hat{A}_{\infty}|(O,0)=1$ and $|\hat{A}_{\infty}|(\cdot,s)\leq4$
on $[-\frac{\alpha}{2},0]$.

On the other hand, let
$\tilde{\varphi}=Q^{-n}\varphi=(4\pi(-s))^{-\frac{n}{2}}e^{-(1+\frac{s}{Q^{2}\varepsilon})
\frac{{d}^{2}_{\tilde{M}}(\bar{x},\cdot)}{4(-s)}-\frac{n}{2\varepsilon}(\bar{t}+Q^{-2}s)},$
note that
\begin{eqnarray*} & &
|H+(1+\frac{1}{\varepsilon}(t-\bar{t}))\frac{d_{\bar{M}}(\bar{x},\cdot)\bar{\nabla}^{\perp}
d_{\bar{M}}(\bar{x},\cdot)}{2(\bar{t}-t)}|_{\bar{g}}^{2}Q^{-2}
=|\tilde{H}-(1+\frac{s}{Q^{2}\varepsilon})\frac{d_{\tilde{M}}(\bar{x},\cdot)\tilde{\nabla}^{\perp}
d_{\tilde{M}}(\bar{x},\cdot) }{2s}|_{\tilde{g}}^{2},\\[2mm]
& &
\psi=(1-\frac{Q^{-2}d_{\tilde{M}}(\bar{x},\cdot)^{2}+3n\bar{t}+3nQ^{-2}s}{\rho^{2}})^{3}_{+}\rightarrow
1,\\[2mm] & & \tilde{\varphi}\rightarrow
(4\pi(-s))^{-\frac{n}{2}}e^{-\frac{|\cdot|^{2}}{4(-s)}}\ \ \ \
\text{and} \ \ \ \varphi\psi dv=\tilde{\varphi}\psi
d\tilde{v}.\end{eqnarray*}

Since $\hat{M}_s\subset \tilde{M}_{s}$, by passing (7.11) to limit,
we have
$$
\int_{-\frac{1}{2}\alpha
}^{0}[\int_{\hat{M}_s^{\infty}}|\hat{H}_{\infty}-\frac{x^{\perp}
 }{2s}|^{2}(4\pi(-s))^{-\frac{n}{2}}e^{-\frac{|x|^{2}}{4(-s)}}]ds=0,
$$
where we denote the limit of $\hat{M}_s$ by $\hat{M}_s^{\infty},$
$\hat{H}_{\infty}$ the mean curvature on the limit.  This implies
$$\hat{H}_{\infty}=\frac{x^{\perp}
 }{2s}\ \ \  \ \ \ \ \ \text{for}\ \  s\in[-\frac{\alpha}{2},0].$$
 The boundedness of the second fundamental form on $\hat{M}_0^{\infty}$  implies
 $x^{\perp}\equiv0$ on $\hat{M}_0^{\infty}.$ Since the second fundamental form and its twice covariant derivative
 of $\hat{M}_s^{\infty}$
 are bounded  for $s\in[-\frac{\alpha}{2},0]$, $\hat{M}_s^{\infty}$ are $C^{4-\gamma}$ submanifolds for any $\gamma>0$.
 Moreover by the higher derivative estimates in Theorem 3.2(in Euclidean space),
 $\hat{M}_0^{\infty}$ is smooth.

 Note $0\in \hat{M}_0^{\infty},$ after a orthogonal transformation, we may assume
 $T_0\hat{M}_0^{\infty}=\{(x_1,x_2,\cdots,x_n,0,\cdots,0)\}$.
 Clearly we still have  the condition $x^{\perp}\equiv 0$ on $\hat{M}_0^{\infty}.$
 We may write $\hat{M}_0^{\infty}$ as a graph (at least locally near $0$ )
 $(x',f_1(x'),\cdots,f_{\bar{n}-n}(x'))$ where
 $x'=(x_1,\cdots,x_n).$ Now $x^{\perp}=(x',f_1(x'),\cdots,f_{\bar{n}-n}(x'))^{\perp}\equiv 0$
 implies $\sum\limits_{p=1}^{n}\frac{\partial f_i}{\partial
 x_{p}}x_{p}=f_{i}(x').$ So $f_i$ is homogenous of degree 1.  Since $Df_i(0)=0,$ we conclude $f_i\equiv
 0.$ Hence we know $\hat{M}_0^{\infty}$ is an $n$-dimensional linear
 subspace $\mathbb{R}^{n}$ of
 $\mathbb{R}^{\bar{n}}.$

 This contradicts $|\hat{A}_{\infty}|(O,0)=1$ and we complete the proof of (7.10).

 Note that $B_{\bar{M}}(\bar{x},\rho)\subseteq B_{\bar{M}}({x}_0,\rho+(2K+1)\varepsilon)\subseteq
 B_{\bar{M}}({x}_0,4\sqrt{\varepsilon})$. Combining  (7.10) and monotonicity formula (7.8), we know
 $$
 \int_{M_{0}\cap
B_{\bar{M}}(x_0,4\sqrt{\varepsilon})}(4\pi\bar{t})^{-\frac{n}{2}}
e^{-(1-\frac{\bar{t}}{\varepsilon})\frac{d^{2}_{\bar{M}}(\bar{x},x)}{4\bar{t}}}dv\geq
\int_{M_{t}}\varphi\psi dv\mid_{t=\bar{t}-\frac{1}{2}\alpha
Q^{-2}}\geq
(1+\beta)e^{-\frac{n}{2\varepsilon}\bar{t}}(1-\frac{3n\bar{t}}{\rho^{2}})^{3}.\eqno{(7.12)}
 $$
By assumption, there is a normal coordinate system  $(y^{1}\cdots
y^{\bar{n}})$ of $\bar{M}$ around $x_0$ with
$T_{x_0}M_0=$span$\{\frac{\partial}{\partial
y^{1}},\cdots,\frac{\partial}{\partial y^{n}}\}$ and a vector
valued function $F:
\{y'=(y^1,\cdots,y^{n})\mid(y^1)^{2}+\cdots+(y^n)^{2}<1\}\rightarrow
\mathbb{R}^{\bar{n}-n} $ with $F(0)=0$, $|DF|(0)=0$,
$|DF|^{2}(y')=\sum_{i,\gamma}\frac{\partial F^{\gamma}}{\partial
y^{i}}\frac{\partial F^{\gamma}}{\partial y^{i}}\leq \delta^{2}$
such that $M_0\cap\{|y'|<1\}=\{(y',F(y'))\mid |y'|<1\}$. Let $P:
\mathbb{R}^{\bar{n}}\rightarrow \mathbb{R}^{{n}}$ be the
orthogonal projection into the first $n$-components. Let
$exp_{x_0}(\bar{y})=\bar{x}$ and $\bar{y}'=P\bar{y}$. For  $x\in
B_{\bar{M}}(x_0,4\sqrt{\varepsilon})$, let $exp_{x_0}({y})={x}$
and ${y}'=P{y}$. Since the curvature of $\bar{M}$ is bounded by
$c_0^{2}$, by comparison theorem  on the ball
$B_{T_{x_0}\bar{M}}(o,4\sqrt{\varepsilon}),$ we have
$$
d_{\bar{M}}(\bar{x},x)\geq
\frac{\sin(4c_0\sqrt{\varepsilon})}{4c_0\sqrt{\varepsilon}}|\bar{y}-y|\geq
(1-3c_0^{2}\varepsilon)|\bar{y}-y|\geq
(1-3c_0^{2}\varepsilon)|\bar{y}'-y'|. \eqno{(7.13)}
$$
On the other hand, also by comparison theorem, the Riemannian volume
element $dv$ of $M_0$ satisfies
$$
\exp_{x_0}^{*}dv\leq [\frac{\sinh (c_0
d_{\bar{M}}({x}_0,\cdot))}{c_0
d_{\bar{M}}({x}_0,\cdot)}]^{n}dv_{\exp_{x_0}^{-1}M_0}\leq
[1+16c_0^{2}\varepsilon]^{n}dv_{\exp_{x_0}^{-1}M_0}\eqno{(7.14)}
$$
whenever $x\in M_0\cap B_{\bar{M}}(x_0,4\sqrt{\varepsilon}).$ By
definition, it is clear that
$$dv_{\exp_{x_0}^{-1}M_0}\leq (1+|DF|^{2})^{\frac{n}{2}}dy^{1}\cdots
dy^{n}\leq (1+\delta^{2})^{\frac{n}{2}}dy^{1}\cdots
dy^{n}.\eqno{(7.15)}$$

Combining (7.13),(7.14) and (7.15), we have
\begin{eqnarray*}
  & & \int_{M_{0}\cap
B_{\bar{M}}(x_0,4\sqrt{\varepsilon})}(4\pi\bar{t})^{-\frac{n}{2}}e^{-(1-\frac{\bar{t}}{\varepsilon})\frac{d^{2}_{\bar{M}}(\bar{x},x)}{4\bar{t}}}dv\\[2mm]
& & \leq (1+\delta^{2})^{\frac{n}{2}}
(1+16c_0^{2}\varepsilon)^{n}(1-\varepsilon)^{-\frac{n}{2}}(1-3c_0^{2}\varepsilon)^{-{n}}\\
& & \ \ \ \ \ \ \
\times\int_{(|y^{1}|^{2}+\cdots+|y^{n}|^{2})^{\frac{1}{2}}<
4\sqrt{\varepsilon}}[\frac{4\pi\bar{t}}{(1-\varepsilon)(1-3c_0^{2}\varepsilon)^{2}}]^{-\frac{n}{2}}
e^{-\frac{|\bar{y}'-y'|^{2}}{\frac{4\bar{t}}{(1-\varepsilon)(1-3c_0^{2}\varepsilon)^{2}}}}dy^{1}\cdots
dy^{n}\\[2mm]
& & \leq (1+\delta^{2})^{\frac{n}{2}}
(1+16c_0^{2}\varepsilon)^{n}(1-\varepsilon)^{-\frac{n}{2}}(1-3c_0^{2}\varepsilon)^{-{n}}.
\end{eqnarray*}
By (7.12) and the fact $\bar{t}\leq \varepsilon^{2}$, we conclude
that
$$
(1+\delta^{2})^{\frac{n}{2}}
(1+16c_0^{2}\varepsilon)^{n}(1-\varepsilon)^{-\frac{n}{2}}(1-3c_0^{2}\varepsilon)^{-{n}}(1-3n
\varepsilon)^{-3}e^{\frac{n\varepsilon}{2}}\geq (1+\beta),
$$
which is a contradiction as $\varepsilon,\delta\rightarrow0$. We
complete the proof of the Theorem.

\end{pf}

\begin{thm} Let $\bar{M}$ be an $\bar{n}$-dimensional manifold
satisfying $\sum\limits_{i=0}^{3}|\bar{\nabla}^{i}\bar{R}m|\leq
c_0^{2}$ and $inj(\bar{M})\geq i_0>0$. Then there is
$\varepsilon>0$ with the following property. Suppose we have a
smooth solution $M_t\subset\bar{M}$ to the MCF  properly embedded
in $B_{\bar{M}}(x_0, r_0)$ for $t\in [0,T]$ where
$r_0<\frac{i_0}{2}$, $0<T\leq \varepsilon^{2} r_0^{2}$. We assume
that at time  zero, $x_0\in M_0,$ and the second fundamental form
satisfies $|A|(x)\leq r_0^{-1}$ on $ M_0\cap B_{\bar{M}}(x_0,r_0)$
and assume $M_0$ is graphic in the ball $B_{\bar{M}}(x_0,r_0).$
Then we have
$$
|A|(x,t)\leq(\varepsilon r_0)^{-1} \eqno{(7.16)}
$$
for any $x\in B_{\bar{M}}(x_0,\varepsilon r_0)\cap M_t$, $t\in
[0,T]$.
\end{thm}
\begin{pf} By scaling we may assume $r_0=1$. By Lemma 7.1, for any
$\delta>0,$ there is $0<r_\delta<1$ such that the connected
component of ${M}_0\cap B_{\bar{M}}(x_0,\frac{1}{96})$ containing
$x_0$ contains a $\delta$-Lipschitz graph of radius $2r_{\delta}$ at
$x_0.$  By our graphic assumption, we conclude that $M_0\cap
B_{\bar{M}}(x_0, r_{\delta})$ is a $\delta$-Lipschitz graph. So
Theorem 7.3 is applicable with radius $r_\delta.$

Consequently, for any $\alpha>0$, there exists an
$\varepsilon_{\alpha}>0$ such that
$$
|A|(x,t)^{2}\leq\frac{\alpha}{t}+\varepsilon_{\alpha}^{-2}
\eqno{(7.17)}
$$
whenever $x\in {M}_t\cap B_{\bar{M}}(x_0,\varepsilon_{\alpha} )$,
$t\in[0,\varepsilon_{\alpha}^{2}]\cap[0,T].$   Let $\alpha$ be a
fixed small constant to be determined later. It turns out that we
only need to choose $\alpha=\alpha(c_0,\bar{n},n)$ finally. Choose
$\varepsilon=\min\{\sqrt{\alpha} \varepsilon_{\alpha}, 10^{-1}\}$.
Then by (7.17) we have
$$
|A|(x,t)^{2}\leq \frac{2\alpha}{t} \eqno{(7.18)}
$$
whenever $x\in {M}_t\cap B_{\bar{M}}(x_0,\varepsilon_{\alpha} )$,
$t\in[0,\varepsilon^{2}]\cap[0,T]$.

\vskip0.3cm \textbf{Claim} \emph{$|A|(x,t)\leq \varepsilon^{-1}$
holds on $M_{t}\cap B_{\bar{M}}(x_0,\varepsilon )$, $t\in
[0,\varepsilon^{2}]\cap[0,T].$} \vskip0.3cm

Suppose $|A|(x_1,t_1)>\varepsilon ^{-1}$ holds for some point
$(x_1,t_1)$, $x_1\in {M}_{t_1}\cap B_{\bar{M}}(x_0,\varepsilon )$,
$t_1\in [0,\varepsilon^{2}]\cap[0,T].$ We can choose another point
$(\bar{x},\bar{t})$,  $\bar{x}\in M_{\bar{t}}\cap
B_{\bar{M}}(x_0,4\varepsilon )$, $\bar{t}\in
[0,\varepsilon^{2}]\cap[0,T]$ such that
$Q=|A|(\bar{x},\bar{t})\geq \varepsilon^{-1}$ and
$$
|A|(x,t)\leq 4 Q \eqno{(7.19)}
$$
whenever $x\in {M}_t$, $d_{\bar{M}}(\bar{x},x)\leq Q^{-1}$, $0\leq
t\leq\bar{t}$.

Actually $(\bar{x},\bar{t})$ can be constructed as the limit of a
finite sequence $(x_i,t_i)$ satisfying  $0\leq t_k\leq t_{k-1}$,
$d_{\bar{M}}(x_0,{x_k})\leq
d_{\bar{M}}(x_0,{x_{k-1}})+|A|({x_{k-1}},t_{k-1})^{-1}$,
$|A|(x_k,t_k)\geq 4|A|(x_{k-1},t_{k-1})$. Since
 $$
|A|(x_k,t_k)\geq 4^{k-1}|A|(x_{1},t_{1})\geq
4^{k-1}\varepsilon^{-1},
 $$
$d_{\bar{M}}(x_0,{x_k})\leq
d_{\bar{M}}(x_0,{x_{1}})+\sum_{i=1}^{\infty}(4^{i-1}|A|(x_{1},t_{1}))^{-1}\leq
3\varepsilon <\frac{1}{2},$ and the solution is smooth, the
sequence must be finite and the last element fits.

Note that $3n\bar{t}Q^{2}\leq 6n\alpha\leq \frac{1}{2}$ by
choosing $\alpha\leq \frac{1}{12n}$. Let
$\psi=(1-\frac{d^{2}_{\bar{M}}(\bar{x},\cdot)+3nt}{Q^{-2}})^{3}_{+}$,
then we have $$(\frac{\partial}{\partial t}-\triangle)\psi \leq
0$$ whenever
$d_{\bar{M}}(\bar{x},\cdot)^{2}<\min\{\frac{1}{c_0^{2}e},i_0^{2}\}$,
$t\in[0,\bar{t}]$. On the other hand, by (3.2), the second
fundamental form satisfies
$$
(\frac{\partial}{\partial t}-\triangle)|A|^{2}\leq -|\nabla
A|^{2}+C(\bar{n})|A|^{4}+C(\bar{n})(1+c_0^{2}) (|A|^{2}+ |A|).
$$
  Hence\begin{equation*}\tag{7.20}
\begin{split}
(\frac{\partial}{\partial t}-\triangle)(\psi |A|^{2})&\leq -|\nabla
A|^{2}\psi+C(\bar{n})|A|^{4}\psi+C(\bar{n})(1+c_0^{2}) (|A|^{2}+
|A|)\psi+4|\nabla A||A||\nabla \psi|\\
&\leq C(\bar{n})|A|^{4}\psi+C(\bar{n})(1+c_0^{2}) (|A|^{2}+
|A|)\psi+4\frac{|\nabla\psi|^{2}}{\psi}|A|^{2}\\
&\leq C(\bar{n})|A|^{4}\psi+C(\bar{n})(1+c_0^{2}) (|A|^{2}+
|A|)\psi+144Q^{2}|A|^{2}\psi^{\frac{1}{3}}
\end{split}
\end{equation*}
on $[0,\bar{t}].$ By (7.19)(7.20), we have
$$
(\frac{\partial}{\partial t}-\triangle)(\psi |A|^{2})\leq
C(\bar{n})Q^{4}+C(\bar{n})(1+c_0^{2})(Q+Q^{2}).
$$
From the maximum principle, it follows
\begin{equation*}
\begin{split}
(\psi |A|^{2})_{\max}\mid_{t=\bar{t}}&\leq 1+
C(\bar{n})Q^{4}\bar{t}+C(\bar{n})(1+c_0^{2})(Q+Q^{2})\bar{t}\\
&\leq 1+2\alpha C(\bar{n})Q^{2}+
C(\bar{n})(1+c_0^{2})(\sqrt{2\alpha\bar{t}}+2\alpha).
\end{split}
\end{equation*}
Note that
$$(\psi |A|^{2})_{\max}\mid_{t=\bar{t}}\geq \psi
|A|^{2}(\bar{x},\bar{t})\geq (1-3nQ^{2}\bar{t})^{3}Q^{2}\geq
(1-18n\alpha)Q^{2},
$$
hence we have
$$
(1-18n\alpha)Q^{2}\leq 1+2\alpha C(\bar{n})Q^{2}+
C(\bar{n})(1+c_0^{2})(\sqrt{2\alpha\bar{t}}+2\alpha).
$$
This implies
$$
Q^{2}\leq \frac{1+
C(\bar{n})(1+c_0^{2})(\sqrt{2\alpha}+2\alpha)}{1-(18n+2C(\bar{n}))\alpha}.
$$
Choosing suitable small $\alpha=\alpha(c_0,\bar{n},n)$, we have
$Q^{2}\leq 2$, which is a contradiction with
$Q^{2}>\varepsilon^{-2}$. So the Claim is proved.

\end{pf}

We remark that in the above theorem the condition that $M_0$ is
graphic in the ball $B_{\bar{M}}(x_0,r_0)$ can be replaced by any
one of the following conditions:

(i) $d_{\bar{g}}(x,y)\geq \frac{d_{g_0}(x,y)}{C}$ for any $x,y\in
M_0\cap B_{\bar{M}}(x_0,r_0);$

(ii)  there is a $\epsilon_0>0$ such that
$B_{\bar{M}}(x_0,\epsilon r_0)\cap M_0$ is connected for any
$\epsilon\leq \epsilon_0.$

\begin{cor} Let $\bar{M}$ be an $\bar{n}$-dimensional complete manifold
satisfying $\sum\limits_{i=0}^{3}|\bar{\nabla}^{i}\bar{R}m|\leq
c_0^{2}$ and $inj(\bar{M})\geq i_0>0$. Let $X_0:M\rightarrow
 \bar{M}$ be an $n$-dimensional isometrically properly embedded  submanifold with bounded
 second fundamental form $|A|\leq c_0$ in $\bar{M}$. We assume $M_0=X_0(M)$ is uniform graphic with some radius $r>0.$
  Suppose $X(x,t)$
 is a smooth solution to the mean curvature flow (1.1)
 on $M\times[0,T_0]$ properly embedded in  $\bar{M}$ with  $X_0$ as initial data.
  Then there is $T_1>0$ depending upon $c_0, i_0, r$ and
  the dimension $\bar{n}$ such that $$|A|(x,t)\leq 2c_0$$ for all $x\in M,$ $0\leq t\leq \min\{T_0,T_1\}$.
\end{cor}
\begin{pf} By Theorem 7.5, there is $\epsilon>0$ such that for any
$x_0\in M$, we have $$ |A|(x,t)\leq \epsilon^{-1}$$  on
$B_{\bar{M}}(x_0,\epsilon)$, $t\in [0,\epsilon^{2}]\cap [0,T].$
Let $[0,\gamma)\subset [0, \epsilon^{2}]\cap[0,T]$  be the maximal
time interval so that the orbit of $x_0$, $x_0^{t}\in
B_{\bar{M}}(x_0,\epsilon)$ for $t\in [0,\gamma].$ Then by the MCF
equation, we know
$$
\frac{d}{dt}d_{\bar{M}}(x_0,x_0^{t})\leq C\epsilon^{-1},
$$
for any $t\in [0,\gamma].$ This implies $\gamma\geq
\frac{\epsilon^{2}}{C}$ for some $C=C(n,\bar{n}).$ Choosing
$\varepsilon=\frac{\epsilon}{\sqrt{C}},$ $T=\min\{T_0,
\varepsilon^{2}\},$ we conclude that the second fundamental forms
 are uniformly bounded by the
constant $\varepsilon^{-1}$ on $M\times[0,T]$. Once the second
fundamental form is bounded, since we assumed
$\sum_{i=0}^{3}|\bar{\nabla}^{i}\bar{R}m|\leq c_0^{2},$ we have
gradient estimate $|\nabla A|\leq \frac{C}{\sqrt{t}},$ and hence
suitable linear growth function with bounded first and second
derivatives can be constructed. Therefore we can apply the maximum
principle to the equation of $|A|$ to conclude a uniform estimate
$|A|\leq 2c_0,$ for any $t\in [0,\frac{1}{C(\bar{n})c_0^{2}}].$ Set
$T_1=\min\{T,\frac{1}{C(\bar{n})c_0^{2}}\}.$ The proof is completed.
\end{pf}
\vskip0.1cm  Theorem 1.3 follows as a corollary of Theorem 1.1 and
Corollary 7.6.

Bing-Long Chen \\
 Department of Mathematics, \\
Sun Yat-Sen university\\
Guangzhou, P.R.China, 510275\\
email: mcscbl@mail.sysu.edu.cn\\

\noindent Le Yin\\
Institute of Mathematical Sciences $\&$ Department of Mathematics,\\
 The Chinese University of Hong Kong,\\
    Hong Kong, P.R.China\\
    email: lyin@math.cuhk.edu.hk
\end{document}